\documentclass[aps,pre,reprint,amsmath,amssymb,longbibliography,floatfix]{revtex4-2}
\usepackage[T1]{fontenc}
\usepackage{lmodern}
\usepackage{mathtools,graphicx,booktabs,enumitem,xcolor,bm}
\usepackage{amsthm}
\usepackage[colorlinks=true,linkcolor=blue!50!black,citecolor=blue!50!black,urlcolor=blue!50!black]{hyperref}
\newtheorem{theorem}{Theorem}
\newtheorem{corollary}[theorem]{Corollary}
\newtheorem{proposition}[theorem]{Proposition}
\newtheorem{lemma}[theorem]{Lemma}
\newtheorem{suppthm}[theorem]{Theorem}
\newtheorem{suppcor}[theorem]{Corollary}
\newtheorem{samplingcor}[theorem]{Corollary}
\newtheorem{suppprop}[theorem]{Proposition}
\newcommand{\C}{\mathbb C}\newcommand{\diag}{\operatorname{diag}}\newcommand{\Prob}{\mathbb P}
\theoremstyle{definition}\newtheorem*{definition}{Definition}
\theoremstyle{remark}\newtheorem*{remark}{Remark}
\newcommand{\R}{\mathbb R}\newcommand{\Fix}{\operatorname{Fix}}
\newcommand{\rank}{\operatorname{rank}}\newcommand{\Diff}{\operatorname{Diff}}
\newcommand{\norm}[1]{\left\lVert#1\right\rVert}
\newcommand{\ip}[2]{\left\langle#1,#2\right\rangle}
\newcommand{\Lip}{\operatorname{Lip}}
\newcommand{\eps}{\varepsilon}\newcommand{\ph}{\varphi}
\begin{document}
\hypersetup{pdftitle={A Takens Embedding under Stochastic Coordinate Erasure: Topological Reconstruction and Finite-Time Prediction},pdfauthor={U Jin Choi}}
\title{A Takens Embedding under Stochastic Coordinate Erasure: Topological Reconstruction and Finite-Time Prediction}
\author{U Jin Choi}
\email[Contact author: ]{ujchoi@kaist.ac.kr}
\affiliation{Department of Mathematical Sciences, Korea Advanced Institute of Science and Technology (KAIST), Daejeon, Republic of Korea}
\date[]{}
\begin{abstract}

The number of surviving measurements does not determine the stability of state reconstruction. We study delay coordinates of a fixed dynamical system when the observed sensor indices can change at each sampling time. For a fixed $C^k$ diffeomorphism ($k\ge2$) of a compact smooth $d$-manifold without boundary and with no periods $1,\ldots,2L-2$, a common open dense set of observation maps yields an embedding for every known length-$L$ mask retaining at least $2d+1$ scalar measurements. The resulting dynamics on each embedded image are conjugate to the original dynamics. This assertion is qualitative: a mask can preserve injectivity while leaving state differences weakly observable. At a specified simple-spectrum periodic point, a sensor--phase rank test is necessary and sufficient for generic local immersion. We also bound retained state separation under specified independent sampling laws. Exact circle examples distinguish the probability of embedding from the probability of well-conditioned inversion. An executable torus example and a simulated nonlinear phase oscillator illustrate how different equal-budget sensor schedules change which forecast leads can be certified. A local dimension estimate can guide window design without certifying it. A whole-block torus control separates dimension estimation from mask-specific state identification. The bounds separate measurement, initialization, and model errors. No claim of exact conjugacy, calibrated coverage, or predictive superiority is made for a fitted model.
\end{abstract}
\maketitle

\section{Introduction}

A sensor network can record the same number of scalar measurements in every window but constrain different state directions from one window to the next. The practical question is not only whether the surviving measurements distinguish states, but also whether small measurement errors force large state errors. We address this question for vector-valued delay observations with known coordinate labels. The analysis proceeds from topological recoverability to observation design, stable inversion, and finite-time prediction.

Takens and subsequent generalized reconstruction results provide the
qualitative framework \cite{takens,sauer,deyle,huke}. In particular,
Deyle and Sugihara's Theorem 7 already fixes the dynamics and treats
prescribed multivariate, nonconsecutive lag patterns in an ambient
prevalence formulation \cite{deyle}. Our qualitative result gives an
intrinsic finite-$C^k$ formulation under explicit short-period exclusion.
Simultaneous validity for a finite family of labeled masks follows by finite
intersection; it is not a new embedding principle. Fixed-system generic
observability and embedding also have precedents
\cite{aeyels1981,gutman}. The conditions and targets differ from
almost-everywhere recovery \cite{baranski} and from the support and
persistence convergence results for missing scalar samples
\cite{dakurah2026}; see Appendix~\ref{app:antecedents}.

Our main mask-specific result, Theorem~\ref{thm:periodic}, gives a
necessary-and-sufficient generic local rank test at a specified simple-spectrum
periodic point. The test accounts for both sensor identity and orbit phase
rather than relying on a scalar count. This is not a global embedding
theorem for systems with short periodic orbits. Stable-delay theory already
distinguishes geometric conditioning from topology \cite{YR,SE}. We use
classical concentration tools to bound retained separation under specified
random mask laws and use deterministic conditioning certificates for realized
masks. An error-budget protocol connects stable inversion to finite-time
prediction without assuming that a data-driven fit inherits conjugacy.

\begin{figure*}[t]
\centering
\includegraphics[width=\textwidth]{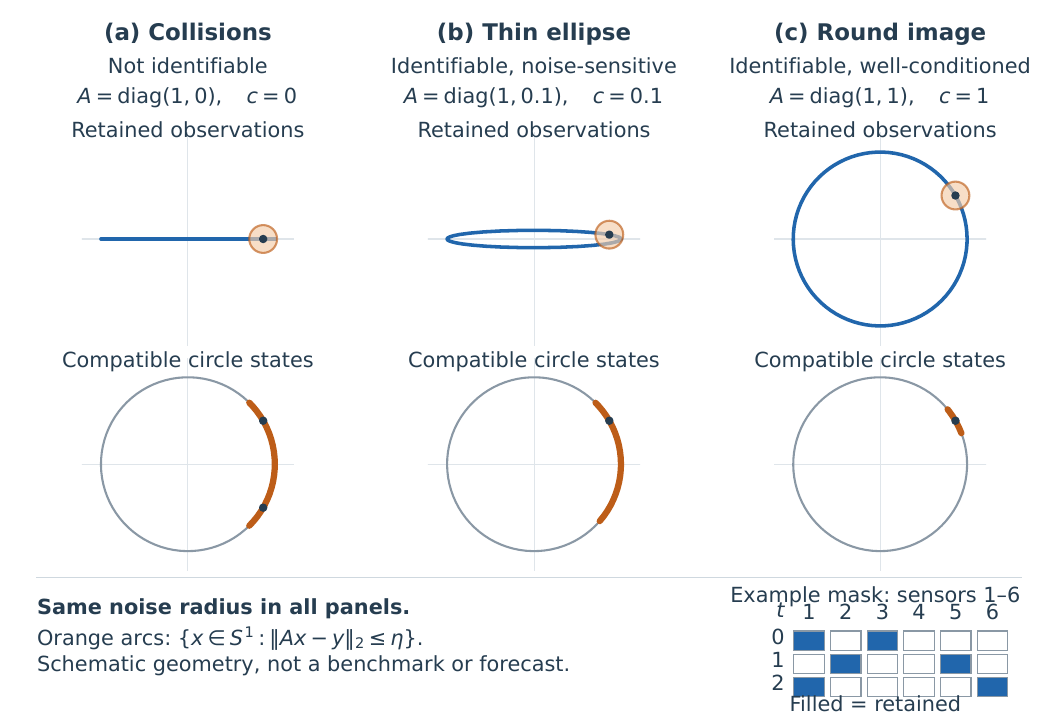}
\caption{Identifiability is not conditioning. These deterministic schematic
maps $x\mapsto Ax$ act on the same unit circle with its Euclidean chordal
metric. Their lower-secant constants are $c=\sigma_{\min}(A)$: (a) a
projection identifies two distinct circle states with the same observation;
(b) a thin ellipse is injective but noise-sensitive; (c) the identity retains
all state separations. Equal-radius observation disks are centered at
$y=Ax_*$, with $x_*=(\sqrt3/2,1/2)$ and illustrative radius $\eta=0.16$.
Orange arcs below are exactly the circle states satisfying
$\|Ax-y\|_2\le\eta$; dots mark $x_*$ and, in (a), its second exact preimage.
The time-by-sensor grid separately illustrates a mask with two retained
coordinates per epoch; it does not generate these three illustrative maps.
In particular, (a) is not an admissible mask of the six-coordinate analytic
example, whose 20 masks are all injective. This is a schematic of physical
observation geometry, not a new benchmark, a learned-decoder guarantee,
or a depiction of dynamical forecast error.}
\label{fig:reader-schematic}
\end{figure*}
\section{Stochastic coordinate-erased delay embedding}
\subsection{Setup and definitions}
Let $M$ be a compact smooth $d$-dimensional manifold without boundary, equipped with a smooth Riemannian metric. Fix $\varphi\in\Diff^k(M)$, $k\ge2$, and $H=(h_1,\ldots,h_D)\in C^k(M,\R^D)$. All genericity statements below use the Whitney $C^k$ topology on the observation space, with $\varphi$ fixed.

At epoch $t$, let $I_t\subset\{1,\ldots,D\}$ be the observed coordinate indices. The observation mask sequence is $\mathbf I=(I_0,\ldots,I_{L-1})$, and the total number of observed scalar entries is
\begin{equation}
 N(\mathbf I)=\sum_{t=0}^{L-1}|I_t|.
\end{equation}
Coordinates within each $I_t$ are ordered increasingly. The per-epoch observation count is $|I_t|$; in the variable-count setting it may vary with $t$ and may be zero.
\begin{definition}[Erased delay map]
For a fixed observation mask sequence $\mathbf I$, define
\begin{equation}\label{eq:delay}
 \Phi_{H,\mathbf I}(x)=\bigl(h_j(\varphi^t x)\bigr)_{0\le t<L,\ j\in I_t}
 \in\R^{N(\mathbf I)}.
\end{equation}
In the fixed-count model, $|I_t|=p=D-r\ge1$ at every epoch and $N(\mathbf I)=Lp$.
\end{definition}
The mask sequence is known to the observer. The theorem concerns each fixed-mask map; it neither discards the coordinate labels nor identifies images associated with different mask sequences.

For an integer $q\ge1$, write $\Fix(\varphi^q)=\{x\in M:\varphi^q(x)=x\}$, where $\varphi^q$ denotes the $q$-fold iterate of $\varphi$. This set consists of periodic points whose minimal period divides $q$; in particular, $\Fix(\varphi)$ is the set of fixed points.

\subsection{Main embedding theorem}
\begin{theorem}[Erased-delay embedding under short-period exclusion]\label{thm:main}
\textbf{(a) Assumptions.} Let $M$, $\varphi\in\Diff^k(M)$, and the observation space $C^k(M,\R^D)$ be as in the setup, with $\varphi$ fixed and $k\ge2$. Fix an integer $L\ge1$ and assume
\begin{equation}\label{eq:aperiodic}
 \Fix(\varphi^q)=\varnothing,\qquad 1\le q\le2L-2.
\end{equation}
For $L=1$, this condition is vacuous. In the fixed-count case, assume $1\le p=D-r<D$ and
\begin{equation}\label{eq:threshold}
 Lp\ge2d+1.
\end{equation}
In the variable-count case, replace the fixed-count requirements by $N(\mathbf I)\ge2d+1$ for each admissible mask sequence.

\textbf{(b) Conclusions.} There is a single open dense subset $\mathcal G\subset C^k(M,\R^D)$ in the Whitney $C^k$ topology such that, for every $H\in\mathcal G$, $\Phi_{H,\mathbf I}$ is a $C^k$ embedding for every admissible mask sequence simultaneously. In particular, this holds for every sequence with $|I_t|=p$ in the fixed-count case. The same $\mathcal G$ can be chosen for the entire finite family
\[
 \begin{split}
 \mathcal A=\bigl\{\mathbf I=(I_0,\ldots,I_{L-1}):\;&I_t\subseteq\{1,\ldots,D\},\\
 &N(\mathbf I)\ge2d+1\bigr\}.
 \end{split}
\]
Thus, after fixing any $H\in\mathcal G$, every probability distribution $\nu$ supported on $\mathcal A$ satisfies
\[
 \Pr_{\mathbf I\sim\nu}\!\left[\Phi_{H,\mathbf I}\text{ is a }C^k\text{ embedding}\right]=1.
\]
The set $\mathcal G$ is independent of $\nu$; independence and uniformity of the masks are not required.
\end{theorem}
\noindent\emph{Proof.} See Appendix~\ref{app:proof-thm-main}.

For each fixed $H\in\mathcal G$ and each admissible mask sequence, embedding means that distinct states have distinct observed delay vectors, $\rank D\Phi_{H,\mathbf I}(x)=d$ for every $x\in M$, and the inverse on $\Phi_{H,\mathbf I}(M)$ is continuous. The probability-one statement concerns only the randomly selected mask sequence: every sequence in the specified support already satisfies this deterministic conclusion. It is not a probability statement about $H$, and it does not guarantee metric preservation, a numerical conditioning bound, or the accuracy of a reconstruction algorithm. The next results distinguish local rank from quantitative inverse sensitivity.

\begin{remark}[Scope of the hypothesis]
Condition \eqref{eq:aperiodic} is sufficient, not necessary. It excludes fixed points whenever $L\ge2$, as well as other short periodic orbits. It is satisfied, for example, by a torus translation $x\mapsto x+\omega$ when $q\omega\notin\mathbb Z^d$ for $1\le q\le2L-2$. Increasing $L$ requires checking the condition over the correspondingly longer range. Genericity of $H$ does not establish the result for an arbitrarily prescribed physical coordinate map.
\end{remark}

\section{Topological conjugacy as the reconstruction target}\label{sec:tc-target}
The primary target is to recover a coordinate representation of the original
dynamics, preserving its topological orbit structure. Numerical accuracy is a
separate question about recovering a particular state and estimating its
subsequent motion. We distinguish an exact, conditional conjugacy construction from
a data-driven estimator, for which no conjugacy is asserted.
Embedding-induced conjugacy is classical \cite{takens,sauer,deyle}; the
contribution here is its mask-aware organization with observation design,
quantitative inversion, and reproducible finite-horizon calculations.

\begin{definition}[Exact dynamical reconstruction]
For a homeomorphism $F:M\to M$, a reconstruction consists of a space $S$, a
homeomorphism $E:M\to S$, and a map $G:S\to S$ such that
\begin{equation}\label{eq:tc-target}
 G E=E F.
\end{equation}
The inverse and intertwining identity are required on the entire declared
state space, not merely on a finite observed trajectory. In the $C^k$ setting,
$E$ and its inverse are $C^k$ maps between manifolds; the reconstructed
manifold need not fill its ambient Euclidean space.
\end{definition}

\begin{theorem}[Conditional exact reconstruction after successful erasure]
\label{thm:tc-core}
Fix the hypotheses of Theorem~\ref{thm:main}, choose $H\in\mathcal G$, and
choose a qualifying mask $\mathbf I$. Put $q=N(\mathbf I)$ and
$E=\Phi_{H,\mathbf I}:M\to S_0=E(M)\subset\mathbb R^q$.
Let $\Psi:S_0\to S=\Psi(S_0)$ be a $C^k$ embedding, where both image
manifolds have dimension $d$. Then $J=\Psi E$ is a $C^k$ diffeomorphism onto
$S$, and the exact reconstructed map
\begin{equation}\label{eq:tc-exact}
 G=J\varphi J^{-1}:S\longrightarrow S
\end{equation}
is a $C^k$ diffeomorphism. In particular,
\begin{equation}\label{eq:tc-iterations}
 GJ=J\varphi,\quad G^nJ=J\varphi^n\ (n\in\mathbb Z),\quad G(S)=S.
\end{equation}
The same conclusion holds if another valid certificate supplies the required
embedding $E$, instead of Theorem~\ref{thm:main}.
\end{theorem}
\begin{proof}
The embedding theorem gives an injective immersion $E$ with a $C^k$ inverse
on its image. The same property holds for $\Psi$, so $J^{-1}=E^{-1}\Psi^{-1}$
is well defined and $C^k$ on $S$. For $x\in M$,
\[
 G(J(x))=J(\varphi(J^{-1}(J(x))))=J(\varphi(x)).
\]
If $G^nJ=J\varphi^n$, then
$G^{n+1}J=GJ\varphi^n=J\varphi^{n+1}$.
Also $G^{-1}=J\varphi^{-1}J^{-1}$, which proves the negative iterates.
For $s=J(x)\in S$, $G(s)=J(\varphi(x))\in S$. Conversely
$s=G(J(\varphi^{-1}(x)))$, so $G(S)=S$. All regularity assertions follow
by composition. No model has been estimated in this argument.
\end{proof}

\paragraph{A constructive formula is not an identification algorithm.}
Equation~\eqref{eq:tc-exact} is executable only when the relevant inverse and
original dynamics are known, or when they have been recovered with additional
proof. Theorem~\ref{thm:main} does not provide those operations from noisy
finite data. In particular, declaring a fitted map to be $G$ does not establish
Eq.~\eqref{eq:tc-target}.

\paragraph{Two legitimate feature choices.}
Taking $\Psi$ to be the identity gives the original erased-delay coordinates.
An ambient diffeomorphism such as $h_a^{(q)}$ in
Eq.~\eqref{eq:ha_definition} restricted to $S_0$ preserves the conclusion.
Alternatively $\Psi(z)=(z,Q(z))$ with a $C^k$ feature vector $Q$ is a graph
embedding: projection onto the first $q$ coordinates is its inverse on the
image, and $D\Psi(z)v=(v,DQ(z)v)$ is injective. If $Q$ is $L_Q$-Lipschitz,
\[
 \|z-z'\|\le\|\Psi(z)-\Psi(z')\|
 \le\sqrt{1+L_Q^2}\|z-z'\|.
\]
This justifies retaining identity coordinates while adding selected quadratic
features. It does not require every scalar basis function to be invertible,
and it cannot separate a pair already collapsed by $E$.

\paragraph{Changing masks.}
For successful embeddings $J_t:M\to S_t$, the exact maps are
\[
 G_t=J_{t+1}\varphi J_t^{-1}:S_t\to S_{t+1},\qquad
 G_tJ_t=J_{t+1}\varphi.
\]
They form a generally nonautonomous coordinate representation, not a single
autonomous map on one fixed image. Multiplying the identities gives
$G_{n-1}\cdots G_0J_0=J_n\varphi^n$. This formulation does not assert
embedding on a failed mask, nor joint topological regularity over a random
mask space without additional hypotheses.

\paragraph{Approximate intertwining and its limited consequence.}
Suppose a fitted map $\widehat G$ preserves $S$, $J^{-1}$ is
$L_{\rm inv}$-Lipschitz from its Euclidean image metric into the
Riemannian distance $d_M$, and
$\sup_{x\in M}\|\widehat G(J(x))-J(\varphi(x))\|\le\epsilon$.
Then, for $\widetilde F=J^{-1}\widehat GJ$,
\[
 d_M(\widetilde F(x),\varphi(x))
 \le L_{\rm inv}\|\widehat G(J(x))-J(\varphi(x))\|
 \le L_{\rm inv}\epsilon.
\]
If $\varphi$ is $L$-Lipschitz in the chosen state metric, this gives
$e_{n+1}\le Le_n+L_{\rm inv}\epsilon$ and hence
$e_n\le L^ne_0+L_{\rm inv}\epsilon\sum_{j=0}^{n-1}L^j$.
This is a conditional finite-time approximation statement, not exact
conjugacy. Even arbitrarily small perturbations need not preserve exact
conjugacy: the identity on a circle and a nonzero small rotation are uniformly
close but cannot be conjugate, since conjugating the identity gives the
identity.

\subsection{Logical scope of topological recoverability}\label{sec:topology-scope}
The following dependency map fixes the meaning of ``recovery'' throughout this paper.
\begin{enumerate}[leftmargin=*]
\item Theorem~\ref{thm:main} and Corollary~\ref{cor:burst} provide global $C^k$ embeddings for a fixed dynamics, a fixed window length, known masks, and a common generic observation law, subject to short-period exclusion. Proposition~\ref{prop:probability} converts the sufficient count event into a lower bound on the probability of this global conclusion; its bound is not an exact embedding probability.
\item Proposition~\ref{prop:periodic}, Theorem~\ref{thm:periodic}, and Corollaries~\ref{new:periodic:2}--\ref{cor:local-reliability} instead concern derivative rank at specified periodic points. Their simple-spectrum assumption does not imply the short-period exclusion or global separation of state pairs.
\item Theorem~\ref{thm:geometry-mask} and its specializations fix the actual full observation map and certify retained separation on a stated sampling event. They require secant information, not genericity. They do not follow from an observation count alone.
\item Once a global embedding is available by either route, Corollary~\ref{cor:classical-conjugacy} gives the classical conjugate dynamics on its image. Proposition~\ref{prop:mask-transition} gives the corresponding transition identities between different successful mask images. Neither construction identifies a fitted regression model with the exact conjugate dynamics.
\end{enumerate}
Thus global topological recoverability, local differential identifiability, quantitative inverse stability, and numerical completion are separate conclusions. The numerical estimation problem is not settled by the preceding guarantees.

\paragraph{What local full rank does imply.}
Let $E:M\to\mathbb R^m$ be $C^k$ and have rank $d$ at $x_0$.
Select $d$ output coordinates whose derivative minor is nonsingular, and
let $P:\mathbb R^m\to\mathbb R^d$ be that coordinate projection.
In a chart at $x_0$, the inverse function theorem gives a neighborhood $U$
on which $P E$ has a $C^k$ inverse. Consequently $E|_U$ is injective, and
its inverse on $E(U)$ is $(P E|_U)^{-1}P$ (with the chart understood).
This proves local embedding, but does not exclude $E(y)=E(x_0)$ for
$y\notin U$. A local image need not be invariant under the dynamics,
so local rank alone does not yield a global conjugate dynamical system.

\subsection{Burst-erasure corollary}
\begin{corollary}[Uniform recovery geometry under burst erasures]\label{cor:burst}
Assume the geometric and regularity hypotheses of Theorem~\ref{thm:main}, including \eqref{eq:aperiodic} for the actual window length $L$.
\begin{enumerate}[label=(\roman*),leftmargin=*]
\item \emph{Fixed per-epoch budget.} If $|I_t|=p$ and $Lp\ge2d+1$, the common open dense set in Theorem~\ref{thm:main} yields embeddings for every mask word, including every word with adversarial coordinate-erasure bursts of length at most $K$. No additional $+K$ window penalty is required by this geometric result.
\item \emph{Only a burst constraint.} Suppose no coordinate is erased for $K+1$ consecutive epochs within the window, with $K\ge0$, and no fixed per-epoch budget is imposed. Then every such word has
\begin{equation}\label{eq:burst-count}
 N(\mathbf I)\ge D\left\lfloor\frac{L}{K+1}\right\rfloor.
\end{equation}
Hence all such words yield embeddings for a common open dense set of $H$ whenever
\begin{equation}\label{eq:burst-bound}
 D\left\lfloor\frac{L}{K+1}\right\rfloor\ge2d+1.
\end{equation}
In particular, a sufficient choice is $L\ge(K+1)\lceil(2d+1)/D\rceil$, provided \eqref{eq:aperiodic} holds for that $L$.
\end{enumerate}
\end{corollary}
\noindent\emph{Proof.} See Appendix~\ref{app:proof-cor-burst}.
\begin{remark}
These are sufficient bounds for noiseless embedding geometry. They are not minimax lower bounds, numerical conditioning estimates, or guarantees for an imputation algorithm. In particular, the additive expression $\lceil(2d+1)/(D-r)\rceil+K$ is not an exact or sharp bound.
\end{remark}

\begin{remark}[Sharpness of the count, not of embedding]
The count in \eqref{eq:burst-count} is attained: observe all $D$ coordinates exactly at epochs $K,2K+1,3K+2,\ldots<L$ and none elsewhere. Every erased run has length at most $K$, and exactly $\lfloor L/(K+1)\rfloor$ epochs are observed. If $L\le K$, the empty word attains zero. Thus the worst-case scalar count is exact; this does not prove necessity or minimax optimality of an embedding or reconstruction threshold.
\end{remark}
\subsection{Periodic-orbit immersion and nonuniform erasures}
\begin{proposition}[Local immersion from consecutive cycle samples]\label{prop:periodic}
This statement is independent of the short-period exclusion in Theorem~\ref{thm:main}. Let $x_0$ have prime period $p_0$, and let $A=D(\varphi^{p_0})_{x_0}$ have $d$ distinct nonzero eigenvalues over $\mathbb C$. Suppose that for some coordinate $j_*$, phase $0\le r_0<p_0$, and integer $q_0\ge0$, this coordinate is observed at every epoch
\begin{equation}\label{eq:consecutive}
 t_i=r_0+(q_0+i)p_0,\qquad i=0,\ldots,d-1,
\end{equation}
with $0\le t_i<L$. Then an open dense set of $H\in C^k(M,\R^D)$ satisfies $\rank D\Phi_{H,\mathbf I}(x_0)=d$.
\end{proposition}
\noindent\emph{Proof.} See Appendix~\ref{app:proof-prop-periodic}.
\begin{remark}
The consecutive cycle indices in \eqref{eq:consecutive} are essential to this proof. Distinct eigenvalues alone do not make a generalized Vandermonde matrix invertible at arbitrary cycle indices: eigenvalues $2,-2$ with indices $0,2$ give identical columns. Proposition~\ref{prop:periodic} gives local immersion at a specified periodic point only. It does not prove global injectivity, cover every mask, or extend Theorem~\ref{thm:main} to systems containing short periodic orbits.
\end{remark}

\subsection{Arbitrary periodic masks and fixed-budget local rank}
\label{sec:periodic}
At a periodic point, repeated measurements need not supply new independent constraints. Which eigenmodes remain distinguishable depends jointly on the sampled cycle indices and the coordinate observed at each phase. The criterion below makes this dependence explicit by grouping measurements with the same coordinate and orbital phase. It uses classical observability and matroid partition \cite{EF,PS} to characterize generic local immersion at one specified periodic point. This local statement neither establishes global injectivity nor removes the short-period hypothesis of Theorem~\ref{thm:main}.

\paragraph{Coordinate identities and aliasing.}
The single-sensor parity obstruction for \(A=\diag(2,-2)\)
is classical \cite[Sec.~V, Eq.~(27)]{dimMo2008}; see also
Ref.~\cite[Sec.~V-A]{PS}. We compare different sensor assignments at
fixed observation epochs and equal total counts. At a two-dimensional
fixed point take \(A=\diag(2,-2)\) and consider epochs \(0\) and \(2\). If both use the same coordinate, the rows are $c_j$ and $4c_j$, so rank is at most one for every $H$. If different coordinates $j,k$ are used, the rows are $c_j$ and $4c_k$, giving rank two generically. Sample times, dynamics, and scalar observation count are identical. The identities of the coordinates change the answer.

Only the two even epochs are included; the intervening odd epoch has no observation. Thus this example is not an exact-positive-budget mask at every original epoch and does not contradict Corollary~\ref{cor:fixed-periodic}.

The linear map $\varphi(x)=\diag(2,-2)x$ on $\R^2$ realizes this local example. No compact global attractor or global embedding claim is involved. In particular, the local threshold $d$ here does not replace the global sufficient threshold $2d+1$.

Let $x_0$ have prime period $p_0$, and suppose
\[
 A=D(\varphi^{p_0})_{x_0}
\]
has $d$ distinct nonzero eigenvalues $\lambda_1,\ldots,\lambda_d$ over $\C$. For each nonempty coordinate--phase group $g=(j,r)$, $0\le r<p_0$, define
\begin{align*}
 Q_g&=\{q\ge0:r+qp_0<L,\ j\in I_{r+qp_0}\},\\
 (V_g)_{qk}&=\lambda_k^q,\qquad q\in Q_g,\quad 1\le k\le d.
\end{align*}
For an eigenmode subset $S\subseteq\{1,\ldots,d\}$ put $r_g(S)=\rank_{\C}V_g[:,S]$.

\begin{suppthm}[Exact generic local-immersion criterion]
\label{thm:periodic}
For a fixed mask and periodic point as above, the following are equivalent:
\begin{enumerate}
\item Some real observation map $H$ has $\rank D\Phi_{H,\mathcal I}(x_0)=d$.
\item An open dense set of $H\in C^k(M,\R^D)$ has this rank, for $k\ge1$.
\item The eigenmode indices can be partitioned into disjoint sets $S_g$ such that each $V_g[:,S_g]$ has full column rank.
\item Every $S\subseteq\{1,\ldots,d\}$ satisfies
\begin{equation}\label{eq:matroid}
 |S|\le\sum_g r_g(S).
\end{equation}
\end{enumerate}
If these conditions fail, every observation map has derivative rank less than $d$ at $x_0$ for this mask.
\end{suppthm}
\noindent\emph{Proof.} See Appendix~\ref{app:new:thm:periodic}.

\begin{suppcor}[Consecutive-cycle sufficient conditions]\label{new:periodic:2}
If one group contains $d$ consecutive cycle exponents $q_0,\ldots,q_0+d-1$, generic local immersion follows. More generally, if every $Q_g$ is consecutive, the criterion is equivalent to $\sum_g|Q_g|\ge d$.
\end{suppcor}
\noindent\emph{Proof.} See Appendix~\ref{app:new:new:periodic:2}.

\begin{suppcor}[Arbitrary switching with a fixed per-epoch budget]
\label{cor:fixed-periodic}
At a specified periodic point with simple nonzero return spectrum, suppose every epoch $t=0,\ldots,L-1$ has exactly $1\le b\le D$ observed coordinates. One common open dense set of observation laws gives local immersion for every such mask if and only if $bL\ge d$. This is a local rank threshold, not a global embedding threshold.
\end{suppcor}
\noindent\emph{Proof.} See Appendix~\ref{app:new:cor:fixed-periodic}.

The same conclusion holds simultaneously at any specified finite collection of periodic points satisfying the spectrum condition. An infinite or positive-dimensional periodic set is not covered by that finite-intersection argument. At least $b$ observations per epoch also suffice by selecting a fixed-budget submask. This corollary removes Proposition~\ref{prop:periodic}'s requirement that a single coordinate recur over $d$ consecutive cycles in the fixed-budget setting; it does not remove Theorem~\ref{thm:main}'s global short-period hypothesis.

\begin{suppcor}[Exact finite-window local reliability]\label{cor:local-reliability}
Fix the periodic point, spectrum, $L,D$, and a finite family $\mathcal F$ of masks. Let $\mathcal G\subseteq\mathcal F$ be the masks satisfying Theorem~\ref{thm:periodic}. There is one open dense observation set $\mathcal U$ such that, for all $H\in\mathcal U$ and any mask distribution $\pi$ on $\mathcal F$,
\begin{equation}\label{eq:localprob}
 \Prob_{\pi}\{\rank D\Phi_{H,\mathcal I}(x_0)=d\}
 =\sum_{\mathcal I\in\mathcal G}\pi(\mathcal I).
\end{equation}
\end{suppcor}
\noindent\emph{Proof.} See Appendix~\ref{app:new:cor:local-reliability}.

For independent coordinate survivals, $\pi(\mathcal I)$ is the product of the appropriate survival and erasure probabilities. For exact-$p$ subset sampling it is the actual subset law, not a product of independent coordinate Bernoulli probabilities.

For the example $A=\diag(2,-2)$, observe one of $D\ge2$ coordinates independently and uniformly at each of $m\ge2$ even epochs. Generic local failure occurs exactly when all $m$ coordinate choices coincide, hence
\[
 \Prob(\text{local rank}<2)=D^{1-m}.
\]
With time-varying independent probabilities $p_{qj}$, the answer is $\sum_{j=1}^D\prod_{q=0}^{m-1}p_{qj}$. These formulas assume a common generic observation law with nonparallel gradients for distinct coordinates. They do not hold for arbitrary prescribed gradients, and are not global embedding probabilities.

\subsection{Count-based probability and its temporal scope}
\begin{proposition}[Independent nonuniform coordinate survival]\label{prop:probability}
Assume \eqref{eq:aperiodic}. Let $s_{t,j}$ be independent Bernoulli variables over all $(t,j)$, with survival probabilities $\rho_{t,j}\in[0,1]$. Set
\[
 N_L=\sum_{t=0}^{L-1}\sum_{j=1}^D s_{t,j},\qquad
 \mu=\mathbb E N_L=\sum_{t,j}\rho_{t,j}.
\]
For every $H$ in the common open dense set for qualifying variable-budget words in Theorem~\ref{thm:main}, if $\mu>2d$, then
\begin{equation}\label{eq:hoeffding}
 \mathbb P(\Phi_{H,\mathbf I}\text{ is not an embedding})
 \begin{aligned}[t]
 &\le\mathbb P(N_L\le2d)\\
 &\le\exp\!\left[-\frac{2(\mu-2d)^2}{LD}\right].
 \end{aligned}
\end{equation}
For time-independent probabilities $\rho_{t,j}=\rho_j$, one has $\mu=L\sum_j\rho_j=LD\bar\rho$.
\end{proposition}
\noindent\emph{Proof.} See Appendix~\ref{app:proof-prop-probability}.
The geometric containment does not require independence; independence is used only for this particular concentration bound. No probability distribution on the observation space is asserted by the word generic.

\begin{suppthm}[Finite horizons versus all time]\label{new:horizon:5}
Suppose a fixed observation law embeds every length-$L$ mask containing at least $2d+1$ scalar entries. If every epoch has exactly $p$ observations and $Lp\ge2d+1$, every window has this geometric property for any known mask sequence, whether random or adversarial.

Under independent Bernoulli sampling, let $q_w$ bound the probability that window $w$ fails the required condition. For any finite collection of $W$ windows,
\[
 \Prob(\text{at least one failure})\le\min\left(1,\sum_{w=1}^W q_w\right).
\]
No independence between overlapping windows is required. If sampling is independent across all epochs and coordinates with constant survival probabilities $\rho_j<1$, and $M$ has at least two points, infinitely many disjoint length-$L$ windows are completely empty almost surely. Consequently an all-time fixed-window embedding guarantee is impossible in that model.
\end{suppthm}
\noindent\emph{Proof.} See Appendix~\ref{app:new:new:horizon:5}.

\begin{remark}
Under an exact-$p$ mask, the within-epoch indicators are dependent; the independent-Bernoulli count model does not apply.
\end{remark}

\section{Observation design and retained separation}\label{sec:design-theory}
\label{sec:mask-stability-extension}

The count condition in the embedding theorem is qualitative. It does not quantify the noise amplification of a realized mask. The following extension fixes the observation law and instead randomizes the observed coordinates. Its ingredients are classical concentration, covering, and conjugacy arguments; the operational contribution is their explicit specialization to coordinate-erased delay maps, including the fixed per-epoch budget with independent masks across epochs. Temporally correlated burst schedules are not covered by that sampling law. Neither a new concentration inequality nor first priority for stable delay embeddings is claimed.

\paragraph{A fixed observation geometry.}
Let $M\subset\mathbb R^a$ contain at least two points and let $F:M\to\mathbb R^N$, $N=LD$, be a fixed full terminal-delay map. For example,
\[
 F(x)=\bigl(h_j(\varphi^{t-L+1}x)\bigr)_{0\le t<L,\,1\le j\le D}.
\]
Assume the certified lower-secant inequality
\begin{equation}\label{eq:full-secant-lower}
 \|F(x)-F(y)\|\ge c_F\|x-y\|\quad (x,y\in M),\qquad c_F>0.
\end{equation}
A fixed $C^1$ embedding of a compact $C^1$ embedded manifold has such a positive constant, but a numerical value requires a separate certificate. No genericity of $H$ or exclusion of short periodic points is needed in the result below once \eqref{eq:full-secant-lower} is established for the actual $F$. If $H(x)=x$ is the physical state observation, the terminal block of $F$ is $x$ itself, so one may take $c_F=1$ immediately. This does not supply the remaining coherence, covering, or reconstruction certificates.

Define the normalized-secant closure and its coordinate coherence by
\begin{equation}\label{eq:secants-coherence}
 \begin{gathered}
 \mathcal S=\overline{\left\{\frac{F(x)-F(y)}{\|F(x)-F(y)\|}:x\ne y\right\}}
 \subset S^{N-1},\\
 \mu=N\sup_{v\in\mathcal S}\|v\|_\infty^2.
\end{gathered}\end{equation}
Thus $1\le\mu\le N$. Let $\mathcal N(\mathcal S,\eta)$ denote the smallest size of a Euclidean $\eta$-net with centers in $\mathcal S$. The closure is compact, so this number is finite. For useful numerical certificates, an explicit upper bound must be provided; compactness alone need not give a useful probability estimate. All geometry and all nets are fixed before sampling the mask.

Here the observation geometry is fixed before the mask is drawn. Randomness selects which entries of its state-separation vectors remain visible; it does not randomize the dynamics or guarantee that a fitted observation map is correct. The lower-secant inequality protects all state pairs, coherence records whether that separating information is concentrated in a few entries, and the covering number measures how many directions must be controlled simultaneously. The resulting probability bound therefore concerns retention of quantitative information, rather than observation count alone. Exact-budget sampling permits dependence within an epoch and assumes independence across epochs.

\begin{suppthm}[Geometry-aware random-mask lower bound]\label{thm:geometry-mask}
Choose one of these two observation laws.
\begin{enumerate}
\item At each epoch $t$, choose a uniform $p$-element subset $I_t\subset\{1,\ldots,D\}$, independently across epochs, where $1\le p\le D$; put $q=p/D$.
\item Retain each of the $N$ coordinates independently with the same probability $q\in(0,1]$.
\end{enumerate}
Let $P_{\mathbf I}$ be the coordinate-restriction map for the resulting mask. For $0<\epsilon<1$, define
\begin{equation}\label{eq:mask-failure-bound}
 B_\epsilon=\min\left\{1,\,
 \mathcal N\left(\mathcal S,\frac{\epsilon\sqrt q}{4}\right)
 \exp\left(-\frac{qN\epsilon^2}{8\mu}\right)\right\}.
\end{equation}
With probability at least $1-B_\epsilon$, simultaneously for all $x,y\in M$,
\begin{equation}\label{eq:mask-stable-lower}
 \|P_{\mathbf I}(F(x)-F(y))\|
 \ge c_F\sqrt{q(1-\epsilon)}\,\|x-y\|.
\end{equation}
In particular, the realized mask is injective. Its inverse-Lipschitz constant is at most
\[ [c_F\sqrt{q(1-\epsilon)}]^{-1}. \] If $F$ is $C^1$ on a manifold, the lower bound also implies immersion, by taking limits along tangent curves. For compact $M$, it is then an embedding.
\end{suppthm}
\noindent\emph{Proof.} See Appendix~\ref{app:new:thm:geometry-mask}.

\begin{suppcor}[Nonuniform survival and variable exact budgets]\label{cor:nonuniform-geometry}
The same proof permits either independent Bernoulli probabilities $\rho_{t,j}\in[0,1]$, or independent uniform exact-$p_t$ sampling per epoch with $\rho_{t,j}=p_t/D$. Define
\[
 \vartheta=\inf_{v\in\mathcal S}\sum_{t,j}\rho_{t,j}v_{t,j}^2.
\]
If $\vartheta>0$, replace every $q$ in \eqref{eq:mask-failure-bound} and \eqref{eq:mask-stable-lower} by $\vartheta$. The resulting conclusion is valid.
\end{suppcor}
\noindent\emph{Proof.} See Appendix~\ref{app:new:cor:nonuniform-geometry}.

\begin{suppcor}[Noise amplification and finitely many rolling windows]\label{cor:noise-rolling}
On the event in Theorem~\ref{thm:geometry-mask}, suppose
\[
 \begin{gathered}y=P_{\mathbf I}F(x)+\xi,\quad\|\xi\|\le\eta_0,\\
 \|P_{\mathbf I}\widetilde F(\widehat x)-y\|\le r,\\
 \sup_{z\in M}\|\widetilde F(z)-F(z)\|\le\beta.\end{gathered}
\]
with $x,\widehat x\in M$. Then
\begin{equation}\label{eq:mask-noise-recovery}
 \|\widehat x-x\|
 \le\frac{\eta_0+r+\beta}{c_F\sqrt{q(1-\epsilon)}}.
\end{equation}
For an exact observation map and a global least-squares minimizer, one may take $r\le\eta_0$.
For $W$ specified rolling windows whose individual mask marginals obey the theorem, the lower-secant and corresponding recovery guarantees hold on all $W$ windows with probability at least
\[
 1-\min\{1,WB_\epsilon\}.
\]
The windows may overlap; independence between windows is unnecessary. This is a finite-horizon guarantee, not an infinite-time assertion.
\end{suppcor}
\noindent\emph{Proof.} See Appendix~\ref{app:new:cor:noise-rolling}.

\paragraph{Interpretation and limitations.}
Equation~\eqref{eq:mask-noise-recovery} separates two sources of difficulty. The numerator is an uncertainty budget: measurement noise, fit residual, and observation-map error. The denominator is a geometric visibility bound for the retained measurements. Increasing the number of measurements helps only insofar as it preserves the state differences that are difficult to observe; high coherence or a large secant covering number can leave the probability estimate uninformative. The guarantee applies to the specified finite collection of windows, even when those windows overlap, but not to an unlimited observation history. Its exact-budget law requires independent mask draws across epochs. A realized correlated burst mask may instead be assessed by a deterministic conditioning certificate. Neither route alone certifies an optimizer or an unknown learned manifold.

For identity observations $H(x)=x$ in physical Euclidean state coordinates, the full terminal map includes the current state as a block, so $c_F\ge1$ follows directly. This does not certify a useful masked lower bound: the required coherence, secant covering bound, and reconstruction residual must still be controlled. Relevant prior includes generalized lags \cite{deyle}, stable delay embedding \cite{YR,SE}, and sampling comparisons and concentration \cite{hoeffding,BM}. The proof in Appendix~\ref{app:new:thm:geometry-mask} is included to make the exact-budget dependence and constants explicit.

\paragraph{A subspace specialization of classical sampling bounds.}
When the full-delay secants lie in a known finite-dimensional subspace,
standard matrix concentration gives a useful alternative to the secant-net bound.
This is a corollary of matrix Chernoff and sampling-without-replacement
comparison, not a new concentration result
\cite{Tropp2011,Tropp2012,GN}.

\begin{samplingcor}[Classical matrix-Chernoff specialization]\label{cor:subspace-certificate}
Let $1\le p\le D$, $N=LD$ and let $U\subset\mathbb R^N$ be a fixed subspace of dimension
$r\ge1$. Choose $Q\in\mathbb R^{N\times r}$ with $Q^\top Q=I_r$ and
$\operatorname{range}(Q)=U$. With rows indexed by $(t,j)$, set
\[
 b=\max_{t,j}\|e_{t,j}^{\top}Q\|_2^2,
 \qquad \mu_U=Nb,\qquad q=p/D.
\]
For each epoch $t$, independently choose a uniform $p$-element subset
$I_t\subset\{1,\ldots,D\}$, and let $P_{\mathbf I}$ restrict to the retained
coordinates. For $0<\epsilon<1$, define
$h_-(\epsilon)=\epsilon+(1-\epsilon)\log(1-\epsilon)$. Then
\[
\begin{split}
&\mathbb P\!\left\{\lambda_{\min}
 (Q^\top P_{\mathbf I}^\top P_{\mathbf I}Q)
 \le q(1-\epsilon)\right\}\\
&\hspace{1em}\le
 \min\left\{1,\;r\exp\!\left[-\frac{qN}{\mu_U}
 h_-(\epsilon)\right]\right\}\\
 &\hspace{1em}\le \min\left\{1,\;r\exp\!\left[-\frac{qN\epsilon^2}{2\mu_U}\right]\right\}.
\end{split}
\]
The same bounds hold under independent Bernoulli retention with probability
$q$. If every secant $F(x)-F(y)$ lies in $U$ and
$\|F(x)-F(y)\|\ge c_F\|x-y\|$, the complementary event implies
\[
 \|P_{\mathbf I}(F(x)-F(y))\|
 \ge c_F\sqrt{q(1-\epsilon)}\,\|x-y\|
.
\]
This holds for every $x,y$.
If $p=D$, the retained map is the full map and the failure probability is
exactly zero.
\end{samplingcor}

\noindent\emph{Proof.} See Appendix~\ref{app:subspace-certificate}.

\paragraph{Scope.}
The printed without-replacement row-sampling theorem
\cite[Thm.~2.2 and Lem.~3.4]{Tropp2011} uses a global uniform subset;
the conditional convex-comparison step above explains its independent
per-epoch extension. No independence of indicators within an epoch is
asserted. Moreover, $\mu_U$ can exceed the secant coherence $\mu$, and a
nonlinear manifold's global secants may span the full ambient space.
This specialization is therefore not uniformly sharper than the general
secant-set result, and tangent-space containment alone is insufficient.

\subsection{A nonvacuous random-mask example}
\label{sec:circle-nonvacuous}
Let $M=S^1\subset\R^2$, $H(x)=x$, and $\varphi(x)=R_\omega x$, with $\omega=2\pi(\sqrt5-1)/2$. Independently at each epoch choose one of the two physical coordinates uniformly. For the terminal full-delay matrix $A_L$ obtained by vertically stacking $R_{-(L-1)\omega},\ldots,I_2$, one has $A_L^\top A_L=LI_2$. Consequently $c_F=\sqrt L$, $\mu=2$, and the normalized secants form an isometric unit circle. A closed net of radius $\epsilon\sqrt q/4$ has
\[\begin{gathered}
 m_\epsilon=\left\lceil\frac{\pi}{2\arcsin(\epsilon/(8\sqrt2))}\right\rceil,
\end{gathered}\]
\[B_\epsilon\le\min\{1,m_\epsilon e^{-L\epsilon^2/16}\}.
\]
Here $D=2,p=1,q=1/2$; no sensor channels are duplicated. At $\epsilon=1/2$, $m_\epsilon=36$. Theorem~\ref{thm:geometry-mask} is genuinely nonvacuous at $L=512$: its failure upper bound is $36e^{-8}<0.01208$, and its all-state retained lower constant is $\sqrt{128}$. This is a fixed-window known-map result, not a learned-model or correlated-burst guarantee.

The same secants lie in a known two-dimensional subspace with $\mu_U=2$. Corollary~\ref{cor:subspace-certificate}, directly Tropp's matrix Chernoff here, gives at $L=64$ the lower constant $4$ and failure bound
\[
 2\exp[-32(1/2+(1/2)\log(1/2))]<0.014751.
\]
Its value is approximately $0.01475021042080273$; the simpler quadratic upper bound is $2e^{-4}<0.036632$. Theorem~\ref{thm:geometry-mask} at $L=64$ still gives $1$, not this sharper probability. This improvement is credited to an existing concentration theorem, not a new principle.

The lower bound prevents flattening in every state-pair direction. Every mask is qualitatively injective for $L\ge2$, yet a direction-dependent mask has squared lower constant at most $L(1/2-1/\pi)<L/4$. Thus the conditioning threshold is genuinely probabilistic. For a global known-map least-squares inverse, total window noise norm $E$ gives state error at most $2E/\sqrt{L(1-\epsilon)/2}$. If instead each retained error is bounded by $\nu$, then $E\le\sqrt L\nu$ and the bound is $2\nu\sqrt{2/(1-\epsilon)}$, with no deterministic $L^{-1/2}$ gain. Appendix~\ref{app:circle-certificate} gives the full geometry, noise normalization, and spectral proof.

\subsection{Analytic example: embedding versus conditioning}\label{sec:analytic-circle}
Write $\Psi_I=P_I F$. Let $M=S^1\subset\mathbb R^2$, let $L=1$, $D=6$, and, for $0<\zeta<1/(2\sqrt3)$, set
\[
 F(u,v)=(u,v,\zeta u,\zeta v,\zeta(u+v),\zeta(u-v)).
\]
The full lower-secant constant is $c_F=\sqrt{1+3\zeta^2}$. Every three-coordinate mask is an embedding: the six row directions are $e_1,e_2,e_1,e_2,e_1+e_2,e_1-e_2$, and no direction occurs three times, so every three chosen rows have rank two. All masks therefore meet the qualitative count $3=2d+1$.

Nevertheless, a three-coordinate mask has lower-secant constant at least $1/2$ if and only if it contains both the first and second coordinates. If it contains them, its squared Gram matrix dominates the identity. If it omits the first coordinate, evaluating on the $e_1$ secant direction bounds its squared lower constant by $3\zeta^2<1/4$; omission of the second coordinate is identical using $e_2$. Every unit direction is a secant direction of the circle. Hence uniform exact-three sampling has
\[
 \begin{gathered}\mathbb P\{c_I\ge1/2\}=\frac{\binom41}{\binom63}=\frac15,\\
 \mathbb P\{\Psi_I\text{ is an embedding}\}=1.\end{gathered}
\]
With independent masks over $W$ epochs, every epoch remains an embedding with probability one, whereas the probability that all $W$ lower constants exceed $1/2$ is exactly $(1/5)^W$. One may take any circle rotation as the underlying dynamics; this example uses the admissible window length $L=1$. Thus random temporal coordinate switching is present even though no delay is required in this example. This is an exact analytic counterexample to equating embedding probability with well-conditioned recovery probability. It is not a benchmark experiment and is not a claim that the preceding sufficient concentration bound is sharp.

\subsection{A design criterion for a finite collection of candidate schedules}
When the observation law is known, let $\mathcal D$ be a specified finite
set of feasible schedules. If $\mathcal W(I)$ lists the windows to be used
under schedule $I\in\mathcal D$, let $\underline c_W$ be a certified
lower-secant bound for window $W$. For noise and fit budget at most $\tau$
uniformly over the windows, the worst certified inversion bound is
$\tau/\min_{W\in\mathcal W(I)}\underline c_W$ when the minimum is positive.
Thus maximizing this minimum over feasible schedules minimizes that
uniform upper bound. For a repeating schedule the window list must
include each phase of the cycle. A zero lower bound means that this
criterion supplies no bound for the corresponding window; it does not
by itself establish noninjectivity. The rule compares certificates, not
empirical recovery errors.

\section{Known-map example: exact torus dynamics}
\label{sec:tc-known}
Consider $b(\theta)=(\cos\theta_1,\sin\theta_1,\cos\theta_2,\sin\theta_2)^\top$
and $R=\operatorname{diag}(R_{\omega_1},R_{\omega_2})$ acting on
$\mathbb T_c^2=S^1\times S^1\subset\mathbb R^4$. In a terminal-state delay
window, an observation matrix $C$ gives row blocks $CR^{t-L+1}$.
The known mask selects a stacked matrix $B\in\mathbb R^{m\times4}$. Suppose it has four independent
columns and $0<c\le\sigma_{\min}(B)$. This numerical check depends on
the actual observations, not just their count.
\[
\begin{gathered}
 S_B=B(\mathbb T_c^2),\qquad E_B(b)=Bb,\\
 G_B(y)=BRB^+y\quad(y\in S_B).
\end{gathered}
\]
Here $B^+=(B^\top B)^{-1}B^\top$.

\begin{proposition}[Known-map torus realization]\label{prop:tc-torus}
The map $E_B:\mathbb T_c^2\to S_B$ is a smooth diffeomorphism onto its
image, with inverse $E_B^{-1}(y)=B^+y$. The restricted map $G_B$ is a smooth
diffeomorphism of $S_B$ and is exactly conjugate to the torus rotation $R$.
\end{proposition}
\begin{proof}
Since $B^+B=I_4$, $Bb=Bb'$ implies $b=b'$, and $B$ is injective on every
tangent space. Its inverse on the image is the restriction of the linear map
$B^+$. Thus $E_B$ is a smooth embedding with the claimed inverse. For any
$b\in\mathbb T_c^2$,
\[
 G_B(E_B(b))=BRB^+Bb=BRb=E_B(Rb).
\]
Since $R$ maps $\mathbb T_c^2$ bijectively to itself,
$G_B(S_B)=S_B$, and the inverse on $S_B$ is $BR^{-1}B^+$.
The identity follows in both directions, proving conjugacy.
\end{proof}
The ambient matrix $BRB^+$ has rank four. If $m>4$, it is singular on
$\mathbb R^m$ but is invertible as a map on $S_B$. This example also shows why an ambient
rank defect alone cannot refute conjugacy on an invariant image.

A simple analytic observation design takes $C=I_4$, $L=2$, and retains
coordinates $\{1,3\}$ at both epochs. Up to row permutation $B$ is block
diagonal, with one block for each circle:
\[
 B_j=\begin{pmatrix}\cos\omega_j&\sin\omega_j\\1&0\end{pmatrix}.
\]
Its two rows have unit norm and inner product $\cos\omega_j$. The
characteristic polynomial of $B_j B_j^\top$ is
$(1-\lambda)^2-\cos^2\omega_j$. Thus
$B_j B_j^\top$ has eigenvalues $1\pm\cos\omega_j$ and
$\sigma_{\min}(B)^2=\min_j(1-|\cos\omega_j|)$.
Choosing $\omega_1=\omega_2=\pi/2$ gives $c=1$. Under the same four-measurement
budget, changing the second mask to $\{2,4\}$ at those same angles repeats
each first-epoch measurement and fails to distinguish states.
This is a deterministic known-map calculation, not an application of the
generic embedding theorem or a learned-model experiment.

\subsection{Stable initialization on the invariant image}
With observations $y=Bb+e$ and $\|e\|\le\eta$, exact linear inversion gives
$v=B^+y$ and $\|v-b\|\le\delta=\eta/c$, where
$c\le\sigma_{\min}(B)$ is certified. Write $v=(v_1,v_2)$ and
$b=(b_1,b_2)$ in two-dimensional blocks. If $\delta<1$, both $v_j$ are
nonzero, and define $p_j=v_j/\|v_j\|$. Then $p\in\mathbb T_c^2$ and
$Bp\in S_B$, so the exact known-map dynamics preserves the reconstructed
image after initialization. The initial estimate need not be the true state.

\begin{proposition}[Product-circle normalization bound]\label{prop:tc-normalization}
Under the preceding assumptions,
\begin{align}
 \|p-b\|&\le Q(\delta):=\sqrt{2-2\sqrt{1-\delta^2}},\label{eq:tc-Q}\\
 d_{\mathbb T}(p,b)&\le\arcsin\delta,\label{eq:tc-angle}
\end{align}
where $d_{\mathbb T}$ is the Euclidean norm of the two principal phase
differences. These are special bounds for exact block normalization, not
for an arbitrary constrained least-squares solver.
\end{proposition}
\begin{proof}
Put $\epsilon_j=\|v_j-b_j\|$, so
$\sum_j\epsilon_j^2\le\delta^2<1$. If $b_j^\top v_j\le0$, then
$\|v_j-b_j\|^2\ge1$, a contradiction. Thus the angle
$\theta_j$ between $p_j$ and $b_j$ belongs to $[0,\pi/2)$.
For $r=\|v_j\|$,
\[
 \epsilon_j^2=(r-\cos\theta_j)^2+\sin^2\theta_j
 \ge\sin^2\theta_j.
\]
Consequently $\theta_j\le\arcsin\epsilon_j$ and
$\|p_j-b_j\|^2\le g(\epsilon_j^2)$, where
$g(s)=2(1-\sqrt{1-s})$. This function is increasing and convex, and
$g(0)=0$. Put $U=\sum_j\epsilon_j^2$. If $U=0$, then $p=b$ and both claims are immediate. If $U>0$, convexity gives
$g(\epsilon_j^2)\le(\epsilon_j^2/U)g(U)$. Summing yields
$\|p-b\|^2\le g(U)\le g(\delta^2)$.
Convexity of $\arcsin$ on $[0,\delta]$ likewise gives
$\arcsin\epsilon_j\le(\epsilon_j/\delta)\arcsin\delta$.
Squaring and summing proves Eq.~\eqref{eq:tc-angle}. If $\delta=0$,
$p=b$ and both bounds follow directly.
\end{proof}
For $0<\epsilon_{\rm tar}<\sqrt2$, solving
$Q(\eta/c)\le\epsilon_{\rm tar}$ gives
\[
 c\ge\frac{\eta}{\epsilon_{\rm tar}
 \sqrt{1-\epsilon_{\rm tar}^2/4}}.
\]

\paragraph{Matrix and solver defects.}
For an approximate solve with $\widetilde B=B+\Delta B$, suppose that
$\|\Delta B\|\le\beta$ and $\sigma_{\min}(\widetilde B)\ge\widetilde c>0$.
With $y=Bb+e$, $\|e\|\le\eta$, define the exact normal-equation defect
of the system actually used by
\[
 \widetilde d=\widetilde B^\top(\widetilde B v-y).
\]
Since $\|b\|=\sqrt2$ and $\widetilde B^+\widetilde B=I_4$,
\begin{align}
 \|\widetilde B^+y-b\|
 &\le \|\widetilde B^+\|\,\|e-\Delta B\,b\|
 \le(\eta+\sqrt2\beta)/\widetilde c,\\
 v-\widetilde B^+y
 &=(\widetilde B^\top\widetilde B)^{-1}\widetilde d.
\end{align}
Hence $\|v-b\|\le(\eta+\sqrt2\beta)/\widetilde c+
\|\widetilde d\|/\widetilde c^2$. The right side can replace $\delta$ in
Eq.~\eqref{eq:tc-Q} if it is less than one. A floating-point certificate
requires validated enclosures for these inputs.

\subsection{Finite-horizon sensitivity}
The state above is already at the last observed epoch. Starting from the
normalized estimate $p$, propagate with the known $R$ or with a rotation
$\widehat R$ for which the vector of angle errors has norm at most $\nu$.
Then $\|\widehat R^n(p-b)\|=\|p-b\|$. In the $j$th circle, the chordal
difference between $\widehat R^n b_j$ and $R^n b_j$ has length
$2|\sin(n\Delta\omega_j/2)|\le n|\Delta\omega_j|$.
Adding these two errors gives
\[
 \|\widehat R^np-R^nb\|\le Q(\delta)+n\nu.
\]
For $\nu>0$ and a target $\epsilon_{\rm tar}\ge Q(\delta)$, a sufficient
horizon is $\lfloor(\epsilon_{\rm tar}-Q(\delta))/\nu\rfloor$.
This bound includes initialization error due to a perturbed observation
matrix through $\delta$; it does not infer model error from a fitted trajectory.

\section{From reconstruction to finite-horizon prediction}\label{sec:forecast}
Consider the following quantitative question: how do errors in the recovered terminal state and learned dynamics propagate after observations stop? The results in this section are conditional on the stated accuracy, stability, and domain assumptions. The symbol $a$ denotes the chosen state-coordinate dimension; it need not equal the intrinsic dimension $d$ or the observation dimension $D$. Full proofs of the formal statements are in Appendix~\ref{app:forecast-proofs}.

Reset the forecast origin to $t=0$. Let $x(t)$ be the true state and $\widehat x(t)$ the forecast. All state errors below use the Euclidean norm in explicitly chosen state coordinates. Define
\[
 E(t)=\norm{\widehat x(t)-x(t)},\qquad E(0)\le\delta_0.
\]
If $\ell$ bounds the one-sided expansion of the true dynamics and $b=\eps+\rho$ bounds the combined dynamics error and numerical residual per unit time, the principal estimate is
\begin{equation}
 E(t)\le e^{\ell t}\delta_0+b\phi_\ell(t).\label{fc:eq:headline}
\end{equation}
Here the accumulation factor is
\[
 \phi_\ell(t)=
 \begin{cases}
 \dfrac{e^{\ell t}-1}{\ell},&\ell\ne0,\\[8pt]
 t,&\ell=0.
 \end{cases}
\]
The reconstruction and propagation factors describe different physical mechanisms. Inversion first converts measurement uncertainty into terminal-state uncertainty through the reciprocal lower-secant constant. The factor $e^{\ell t}$ then transports that initial uncertainty through finite-time dynamical expansion or contraction. The accumulation term accounts for errors introduced after prediction begins: a defect introduced earlier has more time to propagate than one introduced near the evaluation time. Thus better sensing, a more accurate dynamics model, and tighter numerical tolerances address different parts of the error budget.

These inequalities specialize classical differential-inequality, logarithmic-norm, and variation-of-constants arguments \cite{soderlind,jentzen}. Their value here is the explicit connection between a conditioned erased-delay representation and the additional quantities required for a forecast certificate. They are not new propagation principles, and they do not supply data-driven estimates of the required constants.

\paragraph{Three distinctions.}
(i) A fixed diffeomorphism in Theorem~\ref{thm:main} need not be a time-$h$ map of an autonomous ODE; the discrete results do not need this extra assumption. (ii) Completing missing values with observations on both sides of a gap is not a forecast. At the forecast origin, all later observations must be withheld. (iii) A missing-entry RMSE averaged over a record is not a bound on $E(0)$, a uniform vector-field error, or a maximum collocation residual.

\subsection{From an erased-delay embedding to terminal-state error}\label{fc:sec:geometry}
To measure state errors in a fixed norm, realize the compact state manifold $M$ as a $C^1$ embedded submanifold of $\R^a$. Let $\Psi:M\to\R^m$ be a fixed $C^1$ embedding, such as one of the qualifying erased-delay maps in Theorem~\ref{thm:main}. Assume $M$ contains at least two points.

\begin{lemma}[Conditioning of a fixed compact embedding]\label{fc:lem:inverse}
The secant constant
\begin{equation}
 c_\Psi:=\inf_{x\ne y\in M}\frac{\norm{\Psi(x)-\Psi(y)}}{\norm{x-y}}
 \quad\text{satisfies}\quad c_\Psi>0.\label{fc:eq:secant}
\end{equation}
Consequently $\Psi^{-1}$ on $\Psi(M)$ is Lipschitz with constant at most $c_\Psi^{-1}$. For any fixed finite collection of such embeddings, the minimum of their constants is positive. For erased-delay maps this fixes the observation law, window length, and finite mask family; it is not a uniform bound over all generic observation laws or unbounded window lengths.
\end{lemma}
\noindent\emph{Proof.} See Appendix~\ref{app:proof-fc-lem-inverse}.

\begin{corollary}[Noise, fit residual, and observation-map error]\label{fc:cor:recovery}
Let $y=\Psi(x_0)+\xi$, $\norm\xi\le\eta$, and let $\widehat x_0\in M$. Suppose an approximate map $\widetilde\Psi$ satisfies
\[
 \sup_{z\in M}\norm{\widetilde\Psi(z)-\Psi(z)}\le\beta,
 \qquad \norm{\widetilde\Psi(\widehat x_0)-y}\le r_{\rm fit}.
\]
Then
\begin{equation}\norm{\widehat x_0-x_0}\le
 \frac{\eta+r_{\rm fit}+\beta}{c_\Psi}.\label{fc:eq:recovery}
\end{equation}
For an exact map and a global least-squares minimizer on $M$, one may take $r_{\rm fit}\le\eta$, giving $2\eta/c_\Psi$.
\end{corollary}
\noindent\emph{Proof.} See Appendix~\ref{app:proof-fc-cor-recovery}.

The factor $c_\Psi^{-1}$ expresses an inverse problem rather than dynamical instability: states that produce nearly indistinguishable measurements are difficult to separate even before a forecast begins. A small observed residual is informative only relative to this geometric scale and to the errors in the measurement and observation map. In particular, exact agreement with observed coordinates cannot by itself establish accurate recovery of missing coordinates. The bound requires a candidate on the specified state manifold and a certified observation map.

\paragraph{The forecast origin is the right endpoint of the data window.}
For fixed $H$, abbreviate the map in \eqref{eq:delay} by $\Psi_{\mathbf I}=\Phi_{H,\mathbf I}$:
\[
 \Psi_{\mathbf I}(x)=\big(h_j(\ph^t x)\big)_{0\le t<L,\ j\in I_t}.
\]
For this window, let $z=\ph^{L-1}x$ denote the state at its right endpoint. The observation map expressed in that state is
\begin{equation}
 \Psi^{\rm end}_{\mathbf I}(z)
 =\Psi_{\mathbf I}(\ph^{-(L-1)}z)
 =\big(h_j(\ph^{t-L+1}z)\big)_{t,j}.\label{fc:eq:enddelay}
\end{equation}
If $\Psi_{\mathbf I}$ is an embedding, so is $\Psi^{\rm end}_{\mathbf I}$ because $\ph^{-(L-1)}$ is a diffeomorphism. Apply \eqref{fc:eq:recovery} using $c_{\rm end}=c_{\Psi^{\rm end}_{\mathbf I}}$. Alternatively, if $L_{\rm win}=\Lip(\ph^{L-1}|_M)$ in the chosen ambient state norm, then
\begin{equation}
 \begin{aligned}
 c_{\rm end}&\ge \frac{c_{\Psi_{\mathbf I}}}{L_{\rm win}},\\
 E(0)&\le \frac{L_{\rm win}(\eta+r_{\rm fit}+\beta)}{c_{\Psi_{\mathbf I}}}.
 \end{aligned}\label{fc:eq:windowcost}
\end{equation}
Indeed $\norm{z-w}\le L_{\rm win}\norm{\ph^{-(L-1)}z-\ph^{-(L-1)}w}$. One may choose $\delta_0$ equal to the right-hand side of \eqref{fc:eq:windowcost}. Thus inversion at the beginning of a window and forecasting from its end have different conditioning.

\paragraph{Windows ending at an observation.}
For Algorithm~1 we use only windows with $I_{L-1}\ne\varnothing$.
If an input record has trailing empty epochs, let
$m=\max\{t:I_t\ne\varnothing\}$ when this set is nonempty.
Dropping epochs after $m$ leaves the stacked readings and the retained
scalar count unchanged. Use $L=m+1$ in Eq.~\eqref{fc:eq:enddelay} and
evaluate the certificates for the resulting terminal map. Empty epochs
before $m$ remain admissible. A fully erased record has no such endpoint.
This convention does not restrict Theorem~\ref{thm:main}.
If the state at a later, unobserved epoch is the target, that epoch is
a positive forecast lead from $m$ and must pass the corresponding error test.
This is a deterministic convention for realized masks. Conditioning on a
nonempty terminal epoch changes a Bernoulli mask law; an unchanged
independence assumption must not be used for the conditioned coordinates.

\paragraph{Off-manifold estimates.}
If $\widehat x_0\notin M$, \eqref{fc:eq:recovery} cannot be applied directly. One sufficient extension is: choose $\bar x_0\in M$ with $\norm{\widehat x_0-\bar x_0}\le\tau$, let $\Psi$ have an $L_\Psi$-Lipschitz extension to a neighborhood, and assume the approximate-map and residual bounds there. Then
\begin{equation}
 \norm{\widehat x_0-x_0}\le \tau+
 \frac{\eta+r_{\rm fit}+\beta+L_\Psi\tau}{c_\Psi}.\label{fc:eq:offmanifold}
\end{equation}
This follows by comparing $\Psi(\bar x_0)$ first with $\Psi(\widehat x_0)$, then with the data, and adding $\norm{\widehat x_0-\bar x_0}$.

\paragraph{What the embedding theorem does not quantify.}
The finite-mask minimum $c_*>0$ applies to one fixed observation law and one fixed window length. It supplies no uniform positive lower bound over all generic observation laws or all window lengths: replacing $\Psi$ by $a\Psi$ scales $c_\Psi$ by $a>0$. A certified numerical value requires a separate conditioning calculation.

\subsection{Continuous-time propagation}\label{fc:sec:continuous}
\begin{theorem}[Reconstruction and dynamics error propagation]\label{fc:thm:continuous}
Let $x,\widehat x:[0,T]\to\R^a$ be absolutely continuous and satisfy, almost everywhere,
\begin{equation}
 \dot x=f(t,x),\qquad
 \dot{\widehat x}=\widehat f(t,\widehat x)+r_{\rm num}(t).\label{fc:eq:odes}
\end{equation}
Suppose an integrable real function $\ell(t)$ satisfies the one-sided Lipschitz inequality
\begin{equation}
 \ip{u-v}{f(t,u)-f(t,v)}\le\ell(t)\norm{u-v}^2\label{fc:eq:osl}
\end{equation}
on the relevant domain, and nonnegative integrable functions $\eps,\rho$ satisfy
\begin{equation}
 \norm{\widehat f(t,\widehat x(t))-f(t,\widehat x(t))}\le\eps(t),
 \qquad\norm{r_{\rm num}(t)}\le\rho(t).\label{fc:eq:defect}
\end{equation}
Put $A(t)=\int_0^t\ell(s)\,ds$ and $b(t)=\eps(t)+\rho(t)$. If $E(0)\le\delta_0$, then
\begin{equation} E(t)\le e^{A(t)}\delta_0+
 \int_0^t e^{A(t)-A(s)}b(s)\,ds.\label{fc:eq:general}
\end{equation}
This estimate holds for every $t\in[0,T]$.
Only the one-sided inequality for the pairs $(\widehat x(t),x(t))$ is needed; bounds over a whole domain are convenient sufficient conditions.
\end{theorem}
\noindent\emph{Proof.} See Appendix~\ref{app:proof-fc-thm-continuous}.

\begin{corollary}[Embedding-to-forecast inequality]\label{fc:cor:pipeline}
If Corollary~\ref{fc:cor:recovery} holds for the terminal map in \eqref{fc:eq:enddelay}, and Theorem~\ref{fc:thm:continuous} holds in the same state coordinates, then
\begin{equation}
\begin{aligned}
 E(t)\le{}&\frac{e^{A(t)}}{c_{\rm end}}(\eta+r_{\rm fit}+\beta)\\
 &+\int_0^t e^{A(t)-A(s)}\big(\eps(s)+\rho(s)\big)\,ds.
\end{aligned}\label{fc:eq:fullcontinuous}
\end{equation}
For constant upper bounds $\ell(t)\le\ell$ and $\eps(t)+\rho(t)\le b$, this yields \eqref{fc:eq:headline} with $\delta_0=(\eta+r_{\rm fit}+\beta)/c_{\rm end}$.
\end{corollary}
\noindent\emph{Proof.} See Appendix~\ref{app:proof-fc-cor-pipeline}.

The two terms in Eq.~\eqref{fc:eq:fullcontinuous} also explain why improved reconstruction need not yield a proportionate extension of useful prediction. Reducing the terminal inverse error reduces the propagated initial contribution, but leaves the accumulated dynamics and numerical defects unchanged. If the latter already exhaust the allowed error budget, improving the initial reconstruction alone is insufficient. Conversely, a more accurate dynamics model cannot restore state information lost through a poorly conditioned initial observation window. The prediction-time conditions in Appendix~\ref{app:forecast-cases} are therefore joint conditions on observation geometry, model accuracy, and finite-time stability in the same state coordinates.

\paragraph{Computable sufficient expansion bound.}
For $f(t,\cdot)\in C^1$ on a convex region containing the relevant segments, define the Euclidean logarithmic norm
\[
 \mu_2(J)=\lambda_{\max}\big((J+J^\top)/2\big).
\]
If $\sup_z\mu_2(D_xf(t,z))\le\ell(t)$, then \eqref{fc:eq:osl} follows by writing $f(t,u)-f(t,v)=\int_0^1D_xf(t,v+\theta(u-v))(u-v)\,d\theta$ and taking the inner product with $u-v$. An ordinary Lipschitz bound $L(t)\ge0$ also works, but discards possible contraction. One may instead use a one-sided bound for $\widehat f$, by writing the error equation using $\widehat f(t,\widehat x)-\widehat f(t,x)$; then the model defect must be bounded at the true state $x(t)$. A uniform bound on a common region supports either version.

\paragraph{Domain validity and continuation.}
All estimates hold while their domain assumptions hold. A sufficient self-consistency check is available. Suppose $\widehat x$ is given on $[0,T]$, $f$ is continuous in time and locally Lipschitz in state on an open neighborhood of the compact closed $R$-tube around it, and the preceding bounds hold on that tube. If the right-hand side $B(t)$ of \eqref{fc:eq:general} satisfies $B(t)<R$ for every $t\in[0,T]$, then the true trajectory starting with $E(0)<R$ cannot exit the tube. Otherwise the first exit would give $R=E(t_*)\le B(t_*)<R$. Compact containment and the usual continuation argument extend the solution to $T$. Bounds computed only on the training points do not establish this condition.

\paragraph{Sharpness as a worst-case bound.}
On $\R$, take $f(t,x)=\ell(t)x$, $\widehat f(t,x)=\ell(t)x+b(t)$, $r_{\rm num}=0$, $x(0)=0$, and $\widehat x(0)=\delta_0\ge0$. Variation of constants gives equality in \eqref{fc:eq:general}. Thus no smaller universal bound follows from these hypotheses alone. This does not imply that any particular nonlinear forecast attains the bound.

\paragraph{A closed state representation is an additional input.}
An erased-delay embedding defines dynamics on its image. Theorem~\ref{thm:main} does not automatically supply a closed autonomous vector field in arbitrary raw observation coordinates. A raw-coordinate learned predictor therefore requires a separate closure justification, or a genuine state embedding, before it can support a forecast certificate (Appendix~\ref{app:forecast-cases}).

\paragraph{Finite-time growth and changing coordinates.}
A negative asymptotic Lyapunov exponent does not rule out strong transient growth: the nonnormal example in Appendix~\ref{app:forecast-cases} has both eigenvalues $-1$ but finite-time error $\delta_0e^{-t}\sqrt{1+K^2t^2}$. Prediction therefore needs a justified finite-time expansion bound, not an asymptotic exponent alone. Exact mask-dependent coordinate changes cancel between successive transitions (Proposition~\ref{prop:mask-transition}); multiplying their condition numbers at every step would overcount representational distortion. Newly introduced approximation errors do not cancel and must be converted through each receiving mask's inverse conditioning. The appendix retains the horizon cases, discrete and collocation bounds, training-coverage conditions, and a protocol for evaluating future learned models.

\subsection{Intrinsic dimension for observation design}
\label{sec:dimension-design}

The dimension entering the embedding conditions is the geometric dimension
\(d\) of the candidate state manifold, not the rank of the global observation
matrix. A data-derived estimate \(\widehat d\) may guide design, but must not be
substituted for \(d\) in a guarantee. If an independent argument gives an upper
bound \(d_U\ge d\), imposing \(2d_U+1\) instead of \(2d+1\) in a scalar-count
condition is conservative; all the other hypotheses of that condition still
apply. A statistical upper bound, if available, gives this implication only on
its stated coverage event.

A simple reference procedure is available when independent state samples in
fixed ambient dimension are masked completely at random and a complete record
has probability \(p_{\rm full}>0\). The retained complete samples have the
original complete-record law, with an expected count of \(Np_{\rm full}\) among \(N\)
records. One can apply a full-data local-PCA procedure at that reduced sample
size, checking its geometry, noise, and local sample requirements
\cite{dimLim2024}. For a fixed radius, their rule selects the least \(k\) for which the trailing eigenvalues of the centered local covariance contribute at most a fraction \(\eta_{\rm dim}\) of its trace. Complete records may be too scarce for reliable local
neighborhoods: with independent coordinate survivals of probabilities
\(\rho_j\), \(p_{\rm full}=\prod_j\rho_j\).
Repeated measurements may improve coverage if they are registered to the
\emph{same} latent state. In the noiseless case, pooling produces a complete record when every
coordinate is present in at least one repeat;
with noise, the repeat averages still require a precision analysis. Successive
states on a trajectory are not repeats of one state.

An alternative to screen empirically is a split-view procedure: construct
candidate neighborhoods using one view, then estimate first and second
moments using another with recorded co-observation counts
\cite{dimCai2016}. Independent masks alone do not make the selected
neighborhood the true state-space ball. This procedure would require control
of neighborhood error and covariance error relative to a geometric eigengap
before reporting a certified dimension. Even when every coordinate pair has positive probability of joint observation,
the collection of single-epoch masked laws need not identify the dimension.
Appendix~\ref{app:dimension-torus} gives different-dimensional tori with
identical single-epoch masked laws but distinguishable linked two-epoch laws. If no justified upper bound is available, report the
design as conditional on its proposed dimension.

For known invertible dynamics that transports each retained coordinate
to a common calibration epoch, a window that supplies every coordinate
can also yield a complete vector. Conditioning on mask-only coverage
preserves the complete-vector law if masks and the independent
initial states are independent. This argument does not apply to
overlapping windows from a single trajectory, or to a general dynamics
that mixes unknown coordinates. The vector law is that of the complete
observation representation; it has the state dimension only if that
representation preserves dimension, for example as a smooth embedding.

\subsection{Algorithm for reconstruction and finite-time forecasting}
\label{sec:algorithm}

\noindent\textbf{Algorithm 1: masked state reconstruction with an error-budget stopping rule.}
Use a common state metric throughout, here the Euclidean norm on the chosen
state representation. Inputs are the labeled observation window, a specified
admissible state set $M$, a terminal-state observation map or approximation
$\widetilde\Psi_I$, a measurement-noise bound $\eta$, a forecast update
$\widehat F_n$, a finite lead limit $N_{\max}$, and a target error $\Delta$.
The bounds in the steps below are additional inputs or must be computed
and justified for the specified model and domain.

\begin{enumerate}[leftmargin=*,label=\textbf{\arabic*.},itemsep=3pt]
\item \textbf{Design or check the mask.} Use an independently obtained
dimension estimate only to propose candidate windows as in
Sec.~\ref{sec:dimension-design}. Do not count imputed entries toward the
observed budget or use them in the fit residual. If sensor schedules are selectable,
compare feasible schedules by the minimum certified lower-secant bound
over their intended windows, each ending at a nonempty epoch.
For a realized record, discard trailing empty epochs, or return
\emph{not certified} if no measurement remains. Use the final epoch of
this shortened window as the origin and construct the exact or approximate
terminal map for that origin, as in Eq.~\eqref{fc:eq:enddelay}.
Stack the retained values $y$ in the same order as its rows.
Obtain a justified lower-secant bound $\underline c_I>0$ for the true
terminal map and a uniform error bound $\beta$ for $\widetilde\Psi_I$ on $M$.
If either bound is unavailable, return \emph{not certified}; the scalar
observation count is not a substitute.

\item \textbf{Reconstruct the terminal state.} Minimize
$\|\widetilde\Psi_I(z)-y\|$ over $z\in M$, or obtain any feasible candidate
$\widehat x_0\in M$ with a justified residual upper bound $r_{\rm fit}$.
Global optimality is unnecessary if that bound is available. Set
\[
 U_0=(\eta+\beta+r_{\rm fit})/\underline c_I
\]
by Eq.~\eqref{fc:eq:recovery}. If $U_0>\Delta$, do not issue a forecast
at the requested tolerance.

\item \textbf{Propagate and check each lead.} Before releasing lead $n+1$,
form a certified enclosure $K_n$ of
$M\cap\overline B(\widehat x_n,U_n)$, for example the known invariant
domain $M$. The true state lies in this set if the preceding radius holds.
Choose a domain $D_n$ containing both $K_n$ and $\widehat x_n$, on which
the true transition is defined and has a justified Lipschitz bound $q_n$.
If $F_n$ is defined only on $M$, require $\widehat x_n\in M$.
Obtain a model-defect bound $\epsilon_n$ at $\widehat x_n$ (or uniformly
on $D_n$), and compute $\widehat x_{n+1}$ with numerical residual bounded
by $\rho_n$. If the domain or any bound cannot be justified, stop without
releasing that lead. Otherwise update
\[
 U_{n+1}=q_n U_n+\epsilon_n+\rho_n.
\]
Release $\widehat x_{n+1}$ and $U_{n+1}$ only if all needed bounds are
available and $U_{n+1}\le\Delta$. Otherwise stop before releasing that
lead. Starting at $n=0$, repeat up to $N_{\max}$.

\item \textbf{Return the admissible prefix.} Report the terminal estimate
if available, the accepted leads with their error radii, and the stopping
reason: a missing certificate, an error bound above tolerance, or the
requested lead limit. An unknown or excessive bound means that the requested
guarantee is unavailable; it does not prove that the actual error is large.
\end{enumerate}

To see the guarantee, let $E_n$ be the true error at lead $n$. The inverse
bound gives $E_0\le U_0$. If $E_n\le U_n$, then
$x_n\in M\cap\overline B(\widehat x_n,U_n)\subseteq K_n\subseteq D_n$,
and $\widehat x_n\in D_n$ by construction. The discrete propagation
inequality therefore gives $E_{n+1}\le q_n E_n+\epsilon_n+\rho_n\le U_{n+1}$,
since $q_n\ge0$. Thus all released leads meet the tolerance under the
declared assumptions. These are deterministic error radii, not
statistical confidence intervals. The algorithm does not learn the
observation map or the dynamics from incomplete data.

\paragraph{Specialization to the known torus map.}
The inverse step can instead solve the retained linear system in least
squares and normalize each two-coordinate block. For a justified linear
solve error budget $\delta<1$, initialize the forecast radius with
$Q(\delta)$ from Eq.~\eqref{eq:tc-Q}, and use the rotation sensitivity
bound $Q(\delta)+n\nu$ in place of the general recursion for exact
rotation arithmetic. For the quarter-turn schedule retaining
$\{1,3\}$ at both epochs, the known retained matrix has smallest singular
value $1$. With exact inversion arithmetic, a noise bound $\eta=0.01$
therefore gives $\delta=0.01$. The identity
$Q(\delta)^2=2\delta^2/(1+\sqrt{1-\delta^2})$ gives
$0.01<Q(0.01)<0.010001$. For angle-error bound $\nu=0.005$ and target
$\Delta=0.025$, the stopping rule releases leads $1$ and $2$ and
withholds lead $3$: at lead $2$ the upper bound is less than $0.020001$,
while at lead $3$ it is greater than $0.025$.
These are sufficient bounds in a specified example, not observed
errors or benchmark results. The accompanying executable test uses
$(t,j,y_{t,j})=(0,1,3/500),(0,3,1),(1,1,124/125),(1,3,0)$.
For reference terminal state $(1,0,0,1)$, the two nonzero noise entries
are $3/500$ and $-1/125$, giving total Euclidean norm $1/100$.
The signed-coordinate linear solve returns
$(124/125,3/500,0,1)$ before block normalization.
For the forecast, the first circle uses a quarter turn followed by the
rational rotation determined by half-angle tangent $1/400$; its angle
error is $2\arctan(1/400)<1/200$. The second circle uses the exact
quarter turn. Both maps preserve the full product of the two circles.
Given allowance $a=\Delta-n\nu$ with $0\le a<\sqrt2$, the code tests
$Q(\delta)\le a$ by the exact rational criterion
$\delta^2\le a^2(1-a^2/4)$, obtained by squaring
$\sqrt{1-\delta^2}\ge 1-a^2/2$ with both sides nonnegative.
A negative allowance cannot pass, and an allowance of at least $\sqrt2$
always passes for $\delta<1$. This reproduces the two released leads
and rejection at lead 3, and also tests empty records, trailing empty
epochs, and a rank-deficient mask.

\section{A nonlinear oscillator certification study}
\label{sec:adler-certification}
We simulated a running Adler phase oscillator~\cite{adler1946} with a
known observation map. The simulated sensor reports one normalized
Cartesian component (in-phase or quadrature), not the uniform phase
used in the calculation.
This example tests whether the same acquisition budget can yield
different certified forecast horizons when noise and model error are
both nonzero. It is a synthetic study of a phase model, not experimental
hardware data or a learned-model benchmark.

\paragraph{Dynamics and measurement.}
Let $a=1/10$, $w=e^{i\phi}$, and $z=h(w)=(w+a)/(1+aw)=e^{i\theta}$.
For $\dot\phi=\Omega$ the physical phase obeys
\begin{equation}
 \dot\theta=\Omega\frac{1+a^2-2a\cos\theta}{1-a^2}.
 \label{eq:adler-phase}
\end{equation}
Indeed, $d\theta/d\phi=(1-a^2)/|1+aw|^2$ and
$|1+aw|=(1-a^2)/|1-az|$. The drive exceeds the coupling
because $1+a^2>2a$. In Adler form, $\psi=\theta+\pi/2$ gives
$\dot\psi=\Delta\omega-B\sin\psi$, with
$\Delta\omega=\Omega(1+a^2)/(1-a^2)$, $B=2a\Omega/(1-a^2)$, and
$\Delta\omega^2-B^2=\Omega^2$. The time unit is one sampling interval.
Take $\alpha=\Omega=\pi/2+\varepsilon$ and $\varepsilon=\pi\sqrt 2/1000$.
The uniform-coordinate increment $\alpha$ divided by $2\pi$ is
irrational, so the sample map has no periodic points.
One real component of $z$ is recorded at each of
$t=-3,-2,-1,0$. For terminal uniform coordinate $w$, the retained
map $F_I$ has rows $s_{j_t}(h(e^{it\alpha}w))$ where
$s_x(z)=\Re z$ and $s_y(z)=\Im z$. The observations are genuine
nonlinear functions of $w$. The scalar count is $4\ge 2(1)+1$;
we do not infer genericity for these prescribed sensors.

\paragraph{Global lower bounds.}
For $u,v\in S^1$,
\begin{align}
 h(u)-h(v)
 &=\frac{(1-a^2)(u-v)}{(1+au)(1+av)},\\
 |h(u)-h(v)|&\le K|u-v|,\\
 |h(u)-h(v)-(u-v)|&\le d_a|u-v|,
\end{align}
where $K=(1+a)/(1-a)=11/9$ and
$d_a=2a(1+a)/(1-a)^2<0.272$. The last bound follows by expanding
$(1-a^2)-(1+au)(1+av)$ and applying $|u|=|v|=1$.
The nominal linear row at lag $t$ is
$s_{j_t}(e^{it\pi/2}(u-v))$.
Replacing it by the actual nonlinear observation difference changes it
by at most $(d_a+|t|\varepsilon)|u-v|$.
The nominal matrix has signed coordinate rows, with $m_x,m_y$ rows
in the two directions. Thus
\begin{equation}
 \frac{\|F_I(u)-F_I(v)\|}{|u-v|}
 \ge \sqrt{\min(m_x,m_y)}
      -\sqrt{\sum_{t=-3}^{0}(d_a+|t|\varepsilon)^2}.
 \label{eq:adler-mask-lower}
\end{equation}
Here $\varepsilon<0.0045$, so the subtracted term is less than $0.571$.
When $\min(m_x,m_y)$ is $1$ or $2$, respectively, valid lower bounds
$c_I$ in the uniform-phase chord are $0.429$ or $0.843$.
For the all-$x$ mask, an improved bound follows from
\begin{equation}
 \Re h(e^{i\phi})=g(\cos\phi),\qquad
 g(s)=\frac{2a+(1+a^2)s}{1+a^2+2as},
\end{equation}
since $g'(s)=(1-a^2)^2/(1+a^2+2as)^2\ge K^{-2}$ on $[-1,1]$.
Norm comparison before and after detuning gives
$c_I\ge K^{-2}(\sqrt{2}-\varepsilon\sqrt{14})>0.93$.
For a direct physical-chord bound, take $c_I/K$.

The alternating masks $xyxy,yxyx$ have only one nominal direction,
so Eq.~\eqref{eq:adler-mask-lower} alone gives no positive lower bound.
At exact quarter-turn resonance their outputs coincide for
$u,v=1,-1$ and $u,v=i,-i$, respectively. At the detuned angle the
same pair yields an upper bound $K\varepsilon\sqrt{14}<0.021$ on
the best uniform-chord lower constant: each evaluation changes by
at most $K|t|\varepsilon$, and $|u-v|=2$.

Nevertheless the actual alternating masks both embed. Each contains
two $x$ readings separated by two epochs. Before applying $g$, the two
rotated cosine rows have unit norm and inner product $\cos(2\alpha)$.
Their Gram matrix in row coordinates has eigenvalues
$1\pm\cos(2\alpha)$, so the smaller singular value is
$\sqrt{1-|\cos(2\alpha)|}=\sqrt{2}\sin\varepsilon$ at this detuning.
Applying $g$ to each component multiplies every separation by at least
$K^{-2}$ by the mean-value theorem. The full four-row separation can
only be larger than the two-row separation. Hence $0.004$ is a valid
lower constant for both alternating schedules, because
$K^{-2}\sqrt{2}\sin\varepsilon>0.004$, using
$0.0044<\varepsilon<0.0045$ and $\sin z\ge z-z^3/6$ on this interval.
A positive global lower constant also implies immersion by taking
limits along differentiable curves.

\paragraph{Certified forecast and protocol.}
Given observations $y$ and a feasible estimate $\hat w\in S^1$
with residual at most $r$, if
$\|y-F_I(w)\|\le\eta$ and $c_I>0$, then
$|\hat w-w|\le(\eta+r)/c_I$ by the triangle inequality.
Both true and approximate rotations stay on $S^1$. For a model bias
$|\hat\alpha-\alpha|\le\nu$ and lead $n$, the physical chord satisfies
\begin{equation}
 |h(e^{in\hat\alpha}\hat w)-h(e^{in\alpha}w)|
 \le U_n=K\left(\frac{\eta+r}{c_I}+n\nu\right).
 \label{eq:adler-forecast}
\end{equation}
We used exact rotation of the mathematical model, so there is no
time-stepping defect in this bound. The term $n\nu$ follows from
$|e^{in(\hat\alpha-\alpha)}-1|\le n\nu$. We released only the
prefix with $U_n\le0.12$, up to lead 20, with terminal abstention
when the bound was unavailable or above tolerance.

For each of 128 terminal phases drawn independently and uniformly on a
$2^{20}$-point equiangular grid (seed 612021), we used the same signed
noise direction for all schedules. Independent random selection chooses
one of two channels at each epoch. The deterministic geometry-selected
schedule $xxyy$ maximizes the specified certified bound over schedules
using each channel twice, with a fixed tie rule. We also compare $xxxx$,
$yyyy$, and $xyxy$. All have exactly four observations.
We tested noise budgets $\eta=0.001,0.01,0.03$;
the four additive components had independent random signs and
magnitude $0.95\eta/2$. An 8192-point phase grid supplied a feasible
estimate for the known map; its optimality was not used.
Residuals were interval-enclosed, with rational upper bounds rounded
outward using a checked comparison. The forecast deliberately uses the biased uniform-phase increment
$\hat\alpha=\alpha+0.006$, so $\nu=0.006$. Both mathematical
rotations are applied analytically; measured errors in the script
are floating-point diagnostics. All settings for the reported run,
including its original coarse bounds, were fixed before evaluation.

\begin{table*}[t]
\caption{Simulated nonlinear oscillator at the largest noise budget
$\eta=0.03$. All schedules use the same four measurements in count.
Terminal RMSE uses physical chord distance; empirical leads are
counted until the first threshold violation, and certified leads use Eq.~\eqref{eq:adler-forecast}.
Lead $-1$ denotes terminal abstention or, in the empirical column,
terminal error above tolerance. A median of $4.5$ is the median of
integer per-replicate horizons, not a fractional-time prediction.}
\label{tab:adler-study}
\begin{ruledtabular}
\begin{tabular}{lrrrr}
Schedule&Terminal RMSE&Certified starts (of 128)&Median certified lead&Median empirical lead\\
\hline
Fixed $x$&0.01062&128&6&17\\
Fixed $y$&0.00966&128&5&18\\
Alternating $xyxy$&0.25031&0&$-1$&17\\
Independent uniform&0.17595&50&$-1$&18\\
Geometry-selected $xxyy$&0.00987&128&4.5&18\\
\end{tabular}
\end{ruledtabular}
\end{table*}

\begin{figure}[tb]
\centering
\includegraphics[width=\columnwidth]{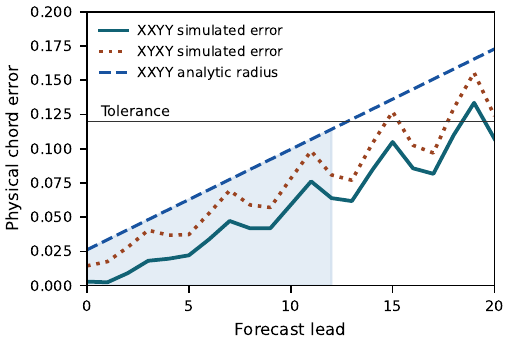}
\caption{One specified replicate (index 0) at $\eta=0.01$.
The shaded area shows the deterministic radius over the accepted
prefix for $xxyy$, not a statistical confidence interval. The dashed
curve continues the bound beyond the release threshold for comparison.
The $xyxy$ error is plotted diagnostically; the original coarse bound
issued no forecasts for that schedule. The refined positive lower
bound also does not meet the specified tolerance in this run.}
\label{fig:adler-trace}
\end{figure}

\paragraph{Empirical outcomes and scope.}
Figure~\ref{fig:adler-trace} illustrates the stopping rule. At $\eta=0.03$, the geometry-selected schedule certified all starts
with leads 4--9, whereas the reported coarse bound released no lead
for the equal-use alternating schedule. Even replacing that bound by
$c_I=0.004$ cannot release the terminal estimate: with the smallest
tested noise budget and any nonnegative residual,
$U_0\ge (11/9)(0.001/0.004)>0.12$. This is failure of that certificate,
not failure of injectivity or proof that the forecast error is large.
At the two larger noise budgets, the upper witness $c_I<0.021$
also gives $U_0>K(0.01)/0.021>0.12$ for any choice of a valid
global lower constant in Eq.~\eqref{eq:adler-forecast}. This does not
rule out a different local uncertainty bound.
The alternating schedule had median terminal
error $0.00825$, but two trials exceeded the tolerance, illustrating why
average performance cannot supply an all-state guarantee.
The selected schedule did not outperform fixed $x$ in certified
horizon, and is compared to alternation under the common two-use
constraint. At $\eta=0.001$ and $0.01$, its median certified leads
were 15 and 12, respectively, with 128 certified starts in each case.
Across the 1920 paired schedule--noise cases, there were 1418 certified
starts and no empirical violations at 17135 released leads, including
terminal estimates. The analytic inequality supplies the guarantee; these
counts only check implementation and show conservatism.
This is one nonchaotic nonlinear system with prescribed maps,
bounded synthetic noise, and known model-error input. It neither
proves a learned identification guarantee nor establishes behavior
under chaotic dynamics or correlated outage laws.

\section{Conclusion}
Known erasure patterns can preserve a delay embedding without preserving a
useful margin for inversion. The periodic rank conditions and the exact
circle and torus examples illustrate how the identities and timing of
measurements can matter as much as their count. A successful embedding
defines a conjugate dynamics on its image, but a data-driven fit does not
automatically recover that map. Finite-time guarantees require a conditioned
terminal reconstruction and separate bounds on model and numerical errors.
The known-model oscillator study shows that equal sensor budgets can
lead to different certificate availability, even when typical empirical
errors overlap. It does not establish that a single schedule is best
across systems. The dimension study adds an independent design input: a local point
estimate is not a replacement for mask-specific state separation.
The exact two-step example in Appendix~\ref{app:dimension-torus} shows
why equal observed counts can give opposite reconstruction decisions.
The protocol separates empirical reconstruction error
from a justified forecast radius and stops when that radius cannot be
validated at the requested tolerance.

\begin{acknowledgments}
The author thanks the Department of Mathematical Sciences at the Korea
Advanced Institute of Science and Technology (KAIST) for providing access
to its research facilities.
\end{acknowledgments}

\section*{Funding}
This research received no research funding.

\section*{AI Use Disclosure}
AI assistance (OpenAI ChatGPT) was used to develop and check mathematical
derivations, prepare and check the example and simulation scripts, and draft and revise
the manuscript. The author set the scope of the analysis and requested
step-by-step derivations and reference checks. The specific model versions
used across sessions were not recorded. AI-assisted checks do not
constitute independent verification. The author is responsible for
verifying the mathematical claims, code, citations, and final text.

\section*{Data Availability}
No external measurements or training datasets were used. The
accompanying source archive includes the analytic-check and torus test
scripts, the whole-block dimension test and its synthetic output,
the phase-oscillator simulation source and seed, and its
per-case and per-lead synthetic results. The package README gives the
commands to regenerate the outputs and the figure.

\appendix
\section{Secondary conditional forecast results}\label{app:forecast-cases}
This appendix continues Section~\ref{sec:forecast}. All statements use the state coordinates, uncertainty budgets, and domain assumptions specified there. Proofs remain in Appendix~\ref{app:forecast-proofs}; the mask-transition proof is in Appendix~\ref{app:new:prop:mask-transition}.

\subsection{Certified prediction time and its dependence on both errors}\label{fc:sec:horizon}
Let $\Delta>0$ be a chosen state-error tolerance and assume $0\le\delta_0\le\Delta$. A \emph{certified horizon} is a time interval on which the bound, and hence the error, remains at most $\Delta$. It is not an upper bound on the actual useful prediction time.

\begin{corollary}[Constant-rate horizons]\label{fc:cor:horizon}
Under the constant bounds of \eqref{fc:eq:headline}, the following horizons are certified, subject throughout to domain validity.
\begin{enumerate}[label=(\roman*)]
\item If $\ell>0$ and $\delta_0+b/\ell>0$,
\begin{equation} T_{\rm cert}=\frac1\ell
 \log\frac{\Delta+b/\ell}{\delta_0+b/\ell}.\label{fc:eq:horizon}
\end{equation}
If $\delta_0=b=0$, the error remains zero.
\item If $\ell=0$, then $E(t)\le\delta_0+bt$. For $b>0$, $T_{\rm cert}=(\Delta-\delta_0)/b$; for $b=0$, every time is certified.
\item If $\ell=-a<0$, put $s=b/a$. Then
\begin{equation}
 E(t)\le s+(\delta_0-s)e^{-at}.\label{fc:eq:contract}
\end{equation}
If $\max\{\delta_0,s\}\le\Delta$, every time is certified. If $\delta_0\le\Delta<s$, then
\begin{equation}
 T_{\rm cert}=\frac1a\log\frac{s-\delta_0}{s-\Delta}.\label{fc:eq:contracthorizon}
\end{equation}
\end{enumerate}
\end{corollary}
\noindent\emph{Proof.} See Appendix~\ref{app:proof-fc-cor-horizon}.

\paragraph{How much does better reconstruction buy?}
For $\ell>0$, a fixed $b$, and $0\le\delta_2\le\delta_1\le\Delta$ with
$\delta_2+b/\ell>0$, reducing the reconstruction bound from $\delta_1$
to $\delta_2$ increases the certificate by
\begin{equation}
 \Delta T_{\rm cert}=\frac1\ell\log
 \frac{\delta_1+b/\ell}{\delta_2+b/\ell}.\label{fc:eq:gain}
\end{equation}
If $b=0$ and $\delta_2=\delta_1/q>0$ with $q>1$, the gain is $(\log q)/\ell$.
If $b=\delta_2=0$, the improved envelope is zero for the entire valid time
interval; if both envelopes are zero there is no finite horizon difference
to report. If $b>0$, even the limit $\delta_2\downarrow0$ gives only
\begin{equation}
 T_{\rm cert}\le\frac1\ell\log\left(1+\frac{\ell\Delta}{b}\right).
 \label{fc:eq:modelceiling}
\end{equation}
This ceiling is for the certificate, not a claim that actual trajectories must separate by that time. It identifies why improving reconstruction without improving the learned dynamics can have diminishing benefit.

\paragraph{A conditional sample-size consequence.}
Suppose, separately from the embedding theorem, one proves $\delta_N\le C_0N^{-\alpha}$ and $b_N\le C_1N^{-\beta}$ on a common valid forecast region, with $\alpha,\beta>0$ and $\ell>0$. Assume that the constants $C_0,C_1,\ell$ are independent of $N$ and that all error and domain hypotheses hold uniformly up to the growing forecast times under consideration, including $t=c\log N$. Put $r=\min\{\alpha,\beta\}$. Then
\begin{equation}
 E_N(t)\le(C_0+C_1/\ell)e^{\ell t}N^{-r}.
\end{equation}
Thus $E_N(c\log N)\to0$ for $0<c<r/\ell$. For a fixed tolerance, a sufficient horizon grows as $(r/\ell)\log N+O(1)$. The proof is direct substitution into \eqref{fc:eq:headline}. These statistical rates are assumptions, not conclusions about a particular estimator. Deteriorating embedding conditioning must be included in $\delta_N$.

\subsection{Finite-time growth, not an asymptotic Lyapunov exponent}
The one-sided bound is useful but may be conservative. A nonlinear variation-of-constants identity exposes the more precise finite-time amplification kernel \cite{jentzen}.

\begin{proposition}[Nonlinear variation of constants]\label{fc:prop:ag}
Let $f$ be autonomous and $C^1$, with flow $\Phi^t$. Suppose all flows and derivatives below exist on a common neighborhood for $0\le s\le t\le T$, including the flows from the straight segment between $x_0$ and $\widehat x_0$. For an absolutely continuous forecast set $d(s)=\dot{\widehat x}(s)-f(\widehat x(s))$. Assume $\norm{d(s)}\le b(s)$ almost everywhere and $\norm{\widehat x_0-x_0}\le\delta_0$. Then
\begin{align}
 \widehat x(t)-\Phi^t(x_0)
 &=\Phi^t(\widehat x_0)-\Phi^t(x_0)\notag\\
 &\quad+\int_0^t D\Phi^{t-s}(\widehat x(s))d(s)\,ds.\label{fc:eq:ag}
\end{align}
If $C_0(t)$ bounds $\norm{D\Phi^t}$ on that initial segment and $C(t,s)\ge\norm{D\Phi^{t-s}(\widehat x(s))}$, then
\begin{equation}
 E(t)\le C_0(t)\delta_0+\int_0^t C(t,s)b(s)\,ds.\label{fc:eq:kernel}
\end{equation}
\end{proposition}
\noindent\emph{Proof.} See Appendix~\ref{app:proof-fc-prop-ag}.

\paragraph{Why $e^{\lambda_{\max}t}$ is not automatically a bound.}
Consider
\[
 A_K=\begin{pmatrix}-1&K\\0&-1\end{pmatrix},\qquad
 e^{tA_K}=e^{-t}\begin{pmatrix}1&Kt\\0&1\end{pmatrix}.
\]
Both eigenvalues, and the largest asymptotic Lyapunov exponent, are $-1$. Yet for initial difference $(0,\delta_0)^\top$,
\begin{equation}
 E(t)=\delta_0e^{-t}\sqrt{1+K^2t^2},
\end{equation}
which can greatly exceed $\delta_0e^{-t}$ on finite intervals. Indeed $\mu_2(A_K)=-1+|K|/2$. A measured long-time Lyapunov exponent cannot replace $\ell$ in \eqref{fc:eq:headline} without a justified finite-time bound, generally including a prefactor and control of the relevant trajectories. Normalizing an empirical prediction time by a Lyapunov time is a reporting convention, not such a proof.

\subsection{Discrete-time forecasts without an ODE assumption}\label{fc:sec:discrete}
\begin{theorem}[Discrete propagation]\label{fc:thm:discrete}
Suppose
\[
 x_{n+1}=F_n(x_n),\qquad
 \widehat x_{n+1}=\widehat F_n(\widehat x_n)+r_n,
\]
and, on the relevant states,
\[
 \begin{gathered}
 \norm{F_n(u)-F_n(v)}\le q_n\norm{u-v},\\
 \norm{\widehat F_n(\widehat x_n)-F_n(\widehat x_n)}\le\eps_n,\\
 \norm{r_n}\le\rho_n,
 \end{gathered}
\]
where $q_n\ge0$. Put $b_n=\eps_n+\rho_n$ and $E_n=\norm{\widehat x_n-x_n}$, and assume $E_0\le\delta_0$. Then
\begin{equation}
 E_n\le\left(\prod_{k=0}^{n-1}q_k\right)\delta_0+
 \sum_{j=0}^{n-1}\left(\prod_{k=j+1}^{n-1}q_k\right)b_j,\label{fc:eq:discrete}
\end{equation}
with empty products equal to one. For constant upper bounds $q_n\le q$ and $b_n\le b$,
\begin{equation}
 E_n\le\begin{cases}
 q^n\delta_0+b(q^n-1)/(q-1),&q\ne1,\\
 \delta_0+nb,&q=1.
 \end{cases}\label{fc:eq:geometric}
\end{equation}
At $q=0$, read the formula for $n\ge1$ as $E_n\le b$ and use $E_0\le\delta_0$ separately.
\end{theorem}
\noindent\emph{Proof.} See Appendix~\ref{app:proof-fc-thm-discrete}.

For $q>1$, $0\le\delta_0\le\Delta$, and $\delta_0+b/(q-1)>0$, the largest integer certified by this envelope is
\begin{equation}
 N_{\rm cert}=\left\lfloor
 \frac{\log[(\Delta+b/(q-1))/(\delta_0+b/(q-1))]}{\log q}
 \right\rfloor.\label{fc:eq:discretehorizon}
\end{equation}
If the physical sampling interval is $h$, the corresponding duration is $hN_{\rm cert}$. For $0\le q<1$ the envelope approaches $b/(1-q)$; if both this value and $\delta_0$ are at most $\Delta$, every step is certified as long as the assumptions remain valid.

\paragraph{Which state is Markov?}
For a fixed erased-delay embedding $\Psi$, the transition on its image is $F_\Psi=\Psi\circ\ph\circ\Psi^{-1}$. Theorem~\ref{thm:main} supports this construction. It does not imply that an arbitrary raw observation $H(x)$ admits a closed autonomous transition $F_H$ or vector field. A quadratic fit in raw standardized observation coordinates would therefore require a separate closure assumption or a forecast formulation in a genuine state embedding. These are not interchangeable without proof. A continuously varying mask would also require a consistent state representation or explicit mask-dependent transitions.

\subsection{Coordinate switching does not itself multiply instability}
For a family $\mathcal A$ of masks, write $\Psi_I=P_I F$ and suppose
\[
 c_I\|x-y\|\le\|\Psi_I(x)-\Psi_I(y)\|\le C_I\|x-y\|,
 \qquad c_I>0.
\]
Let $\varphi:M\to M$ be a map with finite Lipschitz constant $Q$.

\begin{suppprop}[Mask-aware transition identity and stable forcing]\label{prop:mask-transition}
For $I,J\in\mathcal A$, define the exact transition between their embedded images by
\[
 T_{J,I}=\Psi_J\circ\varphi\circ\Psi_I^{-1}.
\]
For every sequence $I_0,\ldots,I_n$,
\begin{equation}\label{eq:switch-telescope}
 T_{I_n,I_{n-1}}\circ\cdots\circ T_{I_1,I_0}
 =\Psi_{I_n}\circ\varphi^n\circ\Psi_{I_0}^{-1},
\end{equation}
and hence
\begin{equation}\label{eq:switch-lipschitz}
 \operatorname{Lip}\bigl(T_{I_n,I_{n-1}}\circ\cdots\circ T_{I_1,I_0}\bigr)
 \le \frac{C_{I_n}}{c_{I_0}}\operatorname{Lip}(\varphi^n)
 \le\frac{C_{I_n}}{c_{I_0}}Q^n.
\end{equation}
Moreover, suppose $\widehat z_k\in\Psi_{I_k}(M)$ and
\[
 \|\widehat z_{k+1}-T_{I_{k+1},I_k}(\widehat z_k)\|\le r_k.
\]
Set $\widehat x_k=\Psi_{I_k}^{-1}(\widehat z_k)$ and $x_k=\varphi^k(x_0)$. Then
\begin{equation}\label{eq:switch-forcing}
 \|\widehat x_n-x_n\|
 \le Q^n\|\widehat x_0-x_0\|
 +\sum_{k=0}^{n-1}Q^{n-1-k}\frac{r_k}{c_{I_{k+1}}}.
\end{equation}
Under the joint good-window event in Corollary~\ref{cor:noise-rolling}, one may use $c_{I_k}\ge c_F\sqrt{q(1-\epsilon)}$ throughout. No additional product of mask condition numbers is required.
\end{suppprop}
\noindent\emph{Proof.} See Appendix~\ref{app:new:prop:mask-transition}.

Changing sensor coordinates changes the representation of the state, not the underlying physical trajectory. For exact mask-dependent embeddings, the intermediate coordinate changes cancel in Eq.~\eqref{eq:switch-telescope}; multiplying a separate coordinate condition number at every step would count representational distortion repeatedly. Approximation errors do not cancel in this way. Each new residual must be converted back to state uncertainty through the receiving mask's inverse conditioning and then propagated by the dynamics, as in Eq.~\eqref{eq:switch-forcing}. This separates genuine dynamical growth from the sensitivity of each measurement representation.

\paragraph{Scope of the transition proposition.}
This is an explicit consequence of coordinate conjugacy, not a new general cocycle theorem. It concerns exact mask-dependent embeddings of the same physical state.
When masks retain different numbers of entries, their ambient spaces
$\mathbb R^{N(I)}$ can have different dimensions. Nevertheless each image
$\Psi_I(M)$ has the same intrinsic dimension as $M$, and $T_{J,I}$ maps
between these image spaces, not between entire ambient Euclidean spaces.
If $\varphi$ is a homeomorphism and both observation maps are embeddings,
$T_{J,I}$ is a homeomorphism with inverse
$\Psi_I\varphi^{-1}\Psi_J^{-1}$. For a changing mask sequence these are
time-dependent coordinate transitions, not one fixed autonomous map. The estimates must be on their corresponding images; arbitrary off-manifold numerical outputs do not meet that assumption. Unknown masks, unknown decoders, approximate transition laws, and learned state representations require additional quantified error bounds. A product of per-step coordinate Lipschitz estimates may be much larger than \eqref{eq:switch-lipschitz}; such artificial growth should not be interpreted as physical instability caused by random sensor switching.

\subsection{A conditional bound for an implicit collocation step}\label{fc:sec:collocation}
Let $g$ be a given increment field. For a pair $(z,z^+)$, define
\begin{equation}
 R_g(z,z^+)=z^+-z-\tfrac12\big[g(z)+g(z^+)\big].\label{fc:eq:trainres}
\end{equation}
The field $g$ need not equal $h$ times a continuous-time generator.
The following error estimate works directly with a discrete transition.

\begin{theorem}[Forecast from a fixed collocation model]\label{fc:thm:collocation}
Let $K\subset\R^a$ be closed, let $F:K\to K$ be the true transition, and let $g:K\to\R^a$ satisfy
\[
 \norm{g(u)-g(v)}\le\kappa\norm{u-v},\qquad 0\le\kappa<2.
\]
Assume $E_0\le\delta_0$, the true and forecast sequences remain in $K$, and the forecast is defined using the fixed field $g$ by
\begin{equation}
\begin{aligned}
 \widehat x_{n+1}&=\widehat x_n+
 \tfrac12\big[g(\widehat x_n)+g(\widehat x_{n+1})\big]+r_n,\\
 \norm{r_n}&\le\rho_n.
\end{aligned}\label{fc:eq:forecastcoll}
\end{equation}
Define the true collocation defect
\begin{equation}
\begin{aligned}
 d_g(x)&=F(x)-x-\tfrac12\big[g(x)+g(F(x))\big],\\
 \sup_{x\in K}\norm{d_g(x)}&\le\eps_{\rm col}.
\end{aligned}\label{fc:eq:coldef}
\end{equation}
Then
\begin{equation}
\begin{aligned}
 E_{n+1}&\le q E_n+\frac{\eps_{\rm col}+\rho_n}{1-\kappa/2},\\
 q&=\frac{1+\kappa/2}{1-\kappa/2}.
\end{aligned}\label{fc:eq:colrec}
\end{equation}
If $\rho_n\le\rho$ and $\kappa>0$, this implies
\begin{equation} E_n\le q^n\delta_0+
 \frac{\eps_{\rm col}+\rho}{\kappa}(q^n-1).\label{fc:eq:colbound}
\end{equation}
If $\kappa=0$, the bound is $\delta_0+n(\eps_{\rm col}+\rho)$.
\end{theorem}
\noindent\emph{Proof.} See Appendix~\ref{app:proof-fc-thm-collocation}.

\paragraph{Existence of the implicit forecast step.}
One sufficient condition for an exact step is that, for each $z\in K$, the map $T_z(w)=z+\tfrac12[g(z)+g(w)]$ maps $K$ into itself. It is a contraction with factor $\kappa/2<1$ on the complete metric space $K$. Iteration $w_{j+1}=T_z(w_j)$ is Cauchy because its successive differences are bounded by a geometric series; its limit is the unique fixed point. Without the self-map condition, $\kappa<2$ alone does not guarantee that a step exists in $K$. For inexact solves, the returned point must remain in $K$ and its actual residual must meet the stated bound. A quadratic polynomial is not globally Lipschitz on all of $\R^a$.

\paragraph{Certificate and conservatism.}
Set $a_{\rm col}=(\eps_{\rm col}+\rho)/\kappa$ when $\kappa>0$. For $0\le\delta_0\le\Delta$ and $\delta_0+a_{\rm col}>0$, the horizon following from \eqref{fc:eq:colbound} is
\begin{equation}
 N_{\rm cert}=\left\lfloor\frac{\log[(\Delta+a_{\rm col})/(\delta_0+a_{\rm col})]}{\log q}\right\rfloor,
\end{equation}
with the zero-error exception as before. This particular proof uses an ordinary Lipschitz constant, so $q\ge1$ even for contracting dynamics. It is a sufficient certificate and may overestimate amplification.

\subsubsection{Training residual, reconstruction error, and coverage}
\begin{proposition}[Conditional defect certificate from finite data]\label{fc:prop:coverage}
Assume additionally that $F$ is $q_F$-Lipschitz on $K$. Let true training inputs $x_i$ form an $r_{\rm cov}$-net of $K$: every $x\in K$ is within $r_{\rm cov}$ of some $x_i$. Suppose reconstructed pairs $\widetilde z_i,\widetilde z_i^+\in K$ satisfy
\[
 \norm{\widetilde z_i-x_i}\le\zeta,\qquad
 \norm{\widetilde z_i^+-F(x_i)}\le\zeta,
\]
and their collocation residuals have maximum norm $R_{\max}$. Then a valid choice in \eqref{fc:eq:coldef} is
\begin{equation}
\begin{aligned}
 \eps_{\rm col}={}&R_{\max}+(2+\kappa)\zeta\\
 &+(1+q_F)(1+\kappa/2)r_{\rm cov}.
\end{aligned}\label{fc:eq:coverage}
\end{equation}
\end{proposition}
\noindent\emph{Proof.} See Appendix~\ref{app:proof-fc-prop-coverage}.

Combining Corollary~\ref{fc:cor:recovery}, Theorem~\ref{fc:thm:collocation}, and Proposition~\ref{fc:prop:coverage}, in one common state representation and for $\kappa>0$, gives the explicit conditional pipeline bound
\begin{equation}\begin{aligned}
 E_n\ \le\;&\frac{q^n}{c_{\rm end}}(\eta+r_{\rm fit}+\beta)\\
 &+\frac{q^n-1}{\kappa}\Big[
 R_{\max}+(2+\kappa)\zeta\\
 &\hspace{28mm}+(1+q_F)(1+\kappa/2)r_{\rm cov}+\rho\Big].
\end{aligned}\label{fc:eq:fullcollocation}
\end{equation}
The first line amplifies terminal reconstruction error. The remaining lines accumulate collocation mismatch, training-state uncertainty, incomplete coverage, and forecast-solver residual. Equation~\eqref{fc:eq:fullcollocation} is a conditional guarantee; applying it to a fitted model requires verifying its constants and domain assumptions.

\paragraph{From mean-square to maximum residual.}
For $m$ pairs $(z_i,z_i^+)$ in $\mathbb R^a$, let
$R_i=R_g(z_i,z_i^+)$ and set
$J_{\rm dyn}=(ma)^{-1}\sum_i\norm{R_i}^2$.
Then $R_{\max}=\max_i\norm{R_i}$ satisfies
\begin{equation}
 R_{\max}\le\sqrt{ma\,J_{\rm dyn}},\label{fc:eq:rmsmax}
\end{equation}
because the maximum squared residual is no greater than the sum of squared residuals. The factor $\sqrt{ma}$ cannot be dropped when converting this average squared residual to a worst-case norm bound. This conversion still provides no bound on $\zeta$, $r_{\rm cov}$, or $q_F$.

\paragraph{Why finite training fit alone is insufficient.}
Even with exact states, a vector-field perturbation can vanish at every finite training point and be arbitrarily large in an unobserved open region, using a smooth bump function supported there. Thus finite residuals alone do not bound the forecast-region error. Smoothness and coverage, as in Proposition~\ref{fc:prop:coverage}, or a different justified generalization theorem are necessary. A single observed trajectory also need not cover a full-dimensional forecast tube.

\subsubsection{When a continuous generator is available}
\begin{corollary}[Vector-field and quadrature components]\label{fc:cor:trapezoid}
Suppose $F$ is the time-$h$ map of $\dot x=f(x)$, true trajectories are $C^3$, and $\norm{x'''(s)}\le M_3$ on each step from $K$. Assume
\[
 \sup_{z\in K}\norm{g(z)-hf(z)}\le h\eps_f.
\]
Then
\begin{equation}
 \eps_{\rm col}\le h\eps_f+\frac{h^3M_3}{12}.\label{fc:eq:quaddefect}
\end{equation}
If $\Lip(g/h)\le L_g$, $0<hL_g<2$, and the forecast step residual is at most $\rho$, one obtains
\begin{equation}
\begin{aligned}
 E_n&\le q^n\delta_0+
 \frac{\eps_f+\rho/h+h^2M_3/12}{L_g}(q^n-1),\\
 q&=\frac{1+hL_g/2}{1-hL_g/2}.
\end{aligned}\label{fc:eq:odecollocation}
\end{equation}
\end{corollary}
\noindent\emph{Proof.} See Appendix~\ref{app:proof-fc-cor-trapezoid}.

For fixed $t=nh$ and $h\downarrow0$,
\[
 \frac{\log q}{h}=L_g+\frac{h^2L_g^3}{12}+O(h^4),\qquad q^n\longrightarrow e^{L_gt}.
\]
Thus the discrete estimate agrees with the continuous bound in the small-step limit when the per-unit-time residual $\rho/h$ is controlled. A fixed nonzero per-step solver tolerance does not give a vanishing accumulated numerical error as $h\downarrow0$.

\subsection{Forecast protocol and interpretation}\label{sec:forecast-protocol}
Theorem~\ref{thm:main} concerns geometric identifiability. The additional results separate two quantitative steps: a conditioned terminal reconstruction gives a bound on $\delta_0$, while dynamics accuracy and finite-time stability propagate it through \eqref{fc:eq:general} or \eqref{fc:eq:colbound}. A small missing-entry average can coexist with an inaccurate terminal state, a biased learned transition, or strong transient growth. The bounds impose no torus-specific geometry; the relevant quantities are conditioning, model defect, and expansion in the forecast region.

For a future autonomous forecast comparison, the following protocol implements the distinction between completion and prediction:
\begin{enumerate}[leftmargin=*,itemsep=4pt]
\item Fix the forecast origin and withhold every later observation before preprocessing, parameter selection, and fitting. Use only a past observation window to estimate the terminal state.
\item Fix the observation scaling and the specified increment field $g$ before forecasting. Use a justified state representation: raw observations require a closed transition, whereas a fixed erased-delay embedding has the transition described in Subsection~\ref{fc:sec:discrete}.
\item Propagate the fixed model by \eqref{fc:eq:forecastcoll}. Record the implicit solver residual at each step. A numerical certificate additionally requires the common-domain and well-posedness conditions of Theorem~\ref{fc:thm:collocation}.
\item Report horizon-dependent error and the first crossing of a prespecified tolerance over common seeds. Report a certified horizon only when the initial error, dynamics defect, amplification constants, and domain conditions have justified bounds.
\end{enumerate}
This procedure uses an implicit collocation step. Offline missing-entry RMSE cannot be substituted for the forecast bound inputs.

\paragraph{Coordinate conversion.}
If standardized coordinates satisfy $x=\mu+Sz$ for a fixed diagonal scaling matrix $S$, then
\[
 \norm{\widehat x-x}\le\norm S_2\norm{\widehat z-z}.
\]
To certify a physical tolerance $\Delta_x$, it is sufficient to use $\Delta_z=\Delta_x/\norm S_2$ in a standardized-state bound. For a Lipschitz output map $H$ with constant $L_H$, output error is at most $L_HE(t)$; a separately estimated output map contributes an additional approximation error. State norms and per-coordinate RMSE tolerances require an explicit conversion.

\section{Relation to prior work and detailed scope}\label{app:antecedents}
A realized mask is a prescribed multivariate, nonconsecutive delay pattern, already within the scope of generalized reconstruction theory \cite{deyle}. Indeed, for each used coordinate set $T_j=\{t:j\in I_t\}$, let $a_j=\min T_j$ and $g_j=h_j\circ\varphi^{a_j}$; then $h_j\circ\varphi^t=g_j\circ\varphi^{t-a_j}$. Their Theorem 7 holds the diffeomorphism fixed and varies the observables by prevalence for prescribed lag relationships in an ambient setting. It is not distinguished from our theorem merely by fixing the dynamics. Theorem~\ref{thm:main} states finite-$C^k$ Whitney-topology genericity on the intrinsic state manifold under a simple short-period exclusion, with the mask quantifiers made explicit. The differing hypotheses should not be identified. Once every mask has an open dense good set in the same observation space, simultaneous validity for a finite mask family follows by intersection. We do not claim a new general embedding principle or first priority for generalized delays.

Stable delay embeddings already distinguish geometric conditioning from topology \cite{YR,SE}. The mask-dependent results specialize classical observability, matroid partition, concentration, and coordinate conjugacy to known erasure masks \cite{EF,PS,hoeffding,BM}. They characterize arbitrary-mask local rank at a specified periodic point and bound noise amplification for a fixed, certified observation geometry. These are additions to this analysis, not claims of new underlying concentration or matroid theorems.

\paragraph{Scalar missing times as a special observation pattern.}
For $D=1$, let $T\subseteq\{0,\ldots,L-1\}$ contain the observed times and choose
\[
 I_t=\begin{cases}\{1\},&t\in T,\\ \varnothing,&t\notin T.\end{cases}
\]
Then $N(\mathbf I)=|T|$ and Eq.~\eqref{eq:delay} becomes
\[
 \Phi_{h,\mathbf I}(x)=\bigl(h(\varphi^t x)\bigr)_{t\in T}.
\]
Thus scalar missing-time windows are included in our variable-count observation model. Under Theorem~\ref{thm:main}'s short-period exclusion, $|T|\ge2d+1$, and genericity of $h$, its conclusion applies to these windows as well as to vector-coordinate masks. For an arithmetic subsequence $T=\{a,a+r,\ldots,a+qr\}$ within the window, set $g=h\circ\varphi^a$ and $F=\varphi^r$. Direct substitution gives
\[
 \Phi_{h,\mathbf I}(x)=\bigl(g(F^i x)\bigr)_{i=0}^{q}.
\]
This identifies a scalar subsequence delay map algebraically; it does not identify a pooled finite-sample estimator with our state-space embedding theorem.

Dakurah--Cisewski-Kehe~\cite[Theorem 4.4 and Corollary 4.5]{dakurah2026} assume nested growing samples, missing counts $l=o(\log m)$, a totally bounded support with a small-ball mass lower bound, and a dependence/avoidance condition on coordinate-disjoint subsamples. Their results control Hausdorff support error and persistence-diagram bottleneck error. Our finite-family, generic-observation statement does not imply these sampling assumptions or their convergence rates. Conversely, convergence of a support or its persistence diagram alone supplies neither an invertible state map nor the intertwining identity in Eq.~\eqref{eq:tc-target}. Hence we regard the present work as an extension of the incomplete-observation research direction to a broader mask model with different conditional guarantees, not a theorem-level generalization of their results.

The closest quantitative antecedents specify different sources of randomness. Yap--Rozell \cite{YR} already give deterministic stable-delay results for a fixed linear observation vector; Eftekhari et al. \cite{SE} randomize coefficients in a prescribed observation-function family. Thus fixed observation geometry is not itself new. Secant covering and noisy inverse bounds are established in manifold embedding \cite{EW}. Tropp's raw row-sampling lemma \cite[Theorem 2.2 and Lemma 3.4]{Tropp2011} already treats uniform exact-size subsets of a fixed incoherent subspace without requiring Hadamard preprocessing. Our per-epoch product-budget law requires the blockwise comparison in Corollary~\ref{cor:subspace-certificate}, using established matrix Chernoff and without-replacement tools \cite{Tropp2012,GN,hoeffding,BM}. Randomized sensor observability \cite{Bopardikar}, intermittent filtering \cite{Sinopoli}, and sampled-coordinate subspace completion \cite{KS} are related but have different sampling laws or target guarantees.

Likewise, validated positive definiteness, singular-value and least-squares bounds \cite{Rump2006,Rump2011,Rump2012}, erasure-robust frames \cite{FM}, and frame-noise amplification \cite{GKK} are established. Our examples use the retained observation geometry explicitly. They establish neither unrestricted all-subset erasure robustness nor a new inverse inequality. Intermittent Kalman filtering connects eigenvalue cycles and irregular sampling to stochastic estimation stability \cite{PS}. Krauss, Lopez, and M\"uller \cite[Sec.~IV A]{krauss2022} explicitly analyze sampled rows $CA^{t_i}$ and relate row dependence to spectral root multiplicities in the irregular single-output linear setting. The repeated-sensor, repeated-phase rows of our local calculation are of this form; our criterion also combines distinct sensor--phase groups and is stated at a specified periodic point. It does not claim priority for the single-output algebra or extend to defective return maps. Ng and Kutz consider uneven sampling in a scalar linear setting; their extension is conditional on a generalized Vandermonde rank conjecture \cite{ngkutz2026}. Ax{\aa}s and Haller also derive local modal Vandermonde structure
for delay embeddings at fixed points, including vector observations
\cite[Sec.~3]{dimAxas2023}. We make no priority claim for this structure. Our rank test is local at a specified periodic point and groups the retained rows by both sensor and orbital phase, with a criterion necessary and sufficient for generic full rank under the stated simple-spectrum assumption.

Missing-value reconstruction is not absent from generalized embedding
research. Carpenter et al.~\cite{carpenter2025gaps} use multiview cross-mapping
to repair ecological time series from dynamically coupled variables,
including short data voids. Their gap-filling method is a relevant
antecedent, not a method excluded by the presence of missingness alone.
Our mask-dependent statements must be compared by their hypotheses and
outputs, not by a blanket claim that earlier work did not handle gaps.

Wu \emph{et al.}~\cite{wu2026} use Gaussian process regression to
learn a full--partial reconstruction map for missing entries. The standard
method uses a complete reference series; a staggered-segment lag variant
also handles synchronized gaps without a fully observed channel.
These are learned imputation results under their data and training setup,
not per-mask embedding or deterministic finite-lead certificates. The
present analysis is not a theorem-level generalization of their results.

Coordinate dropouts arise in distributed sensors, tracking, and telemetry. Interpolation, multivariate imputation, and matrix completion provide useful computational baselines \cite{mice,candes}. Their empirical performance depends on the observation regime and the data geometry; no universal failure of these methods is asserted here.

Partial observations also motivate data-driven prediction in nonlinear physics. Bhat and Munch connect delay representations to recurrent neural networks and approximation error \cite{bhat}. Young and Graham learn delay-coordinate dynamics and state reconstruction from partial observations \cite{young}. Shi et al. use a reservoir model with an explicit time-interval input to forecast trajectories with missing time samples \cite{shi}. These studies establish relevant reconstruction and prediction approaches. Here we study known, time-varying coordinate masks and a simultaneous-mask embedding guarantee under short-period exclusion. Missing time samples, fixed partial observations, and changing coordinate masks are distinct observation models; their empirical comparisons require matching the available information and training protocol.

The mathematical issue is that counting observed scalars is not the same as counting independent constraints. Repeated orbit points can cause repeated observations to contain dependent information. A volume-zero periodic set does not resolve this obstruction to a global embedding. Our sufficient hypothesis excludes short periodic orbits and allows the independence argument to be completed. This is a restricted theorem for a fixed $\varphi$ and generic $H$, not a claim for all dynamical systems or a generic-pair theorem.

\section{Proofs of the embedding and erasure results}\label{app:embedding-proofs}

\subsection{Proof of Theorem~\ref{thm:main}}\label{app:proof-thm-main}
\begin{proof}
Fix one qualifying mask word and write $m=N(\mathbf I)\ge2d+1$.

\paragraph{Step 1: immersion.}
For $d=0$ immersion is automatic; the following tangent-bundle argument is needed only for $d\ge1$. The points $z_t=\varphi^t x$, $0\le t<L$, are distinct. Otherwise $\varphi^{b-a}x=x$ for some $0\le a<b<L$, contradicting \eqref{eq:aperiodic}. For $L=1$ distinctness is automatic. Let
\[
 SM=\{(x,v)\in TM:\|v\|=1\},\qquad\dim SM=2d-1,
\]
and define $F(H,x,v)=D\Phi_{H,\mathbf I}(x)v$. Its components are
\begin{equation}\label{eq:derivative}
 F_{t,j}(H,x,v)=Dh_j(z_t)w_t,\qquad w_t=D(\varphi^t)_xv\ne0.
\end{equation}
For arbitrary numbers $a_{t,j}$, choose covectors $\ell_t$ with $\ell_t(w_t)=1$. Local coordinate functions multiplied by bump functions in disjoint neighborhoods of the $z_t$ realize
\[
 D(\delta h_j)(z_t)=a_{t,j}\ell_t,\qquad j\in I_t.
\]
At one point, the different coordinate functions may be changed independently. Hence
\begin{equation}
 D_HF(H,x,v)[\delta H]=(a_{t,j})_{t,j}.
\end{equation}
The parameter derivative is surjective. Step 3 below uses these perturbations
in finite-dimensional parameter spaces to obtain an open dense set of
observations for which $F_H$ has no zero, hence $\Phi_{H,\mathbf I}$ is
an immersion. It does not apply finite-dimensional transversality
directly to the space of observation laws.

\paragraph{Step 2: secant constraints and overlapping windows.}
On $M^{(2)}=(M\times M)\setminus\Delta$, define
\begin{equation}
 E(H,x,y)=\Phi_{H,\mathbf I}(x)-\Phi_{H,\mathbf I}(y).
\end{equation}
We prove surjectivity of $D_HE$ at every $x\ne y$. If the two orbit windows are disjoint, their function values are freely adjustable by bump functions, and surjectivity follows immediately.

Otherwise $\varphi^a x=\varphi^b y$ for some $0\le a,b<L$, so $y=\varphi^s x$ with $s=a-b$ and $0<|s|\le L-1$. All orbit points in the two windows are indexed by
\[
 V=\{0,\ldots,L-1\}\cup\{s,\ldots,s+L-1\}.
\]
The index span is $L-1+|s|\le2L-2$. Distinct indices in $V$ represent distinct points by \eqref{eq:aperiodic}, including when some indices are negative, since $\varphi$ is invertible.

For each fixed coordinate $j$, put an oriented edge from $t$ to $t+s$ whenever $j\in I_t$. This graph is a subgraph of disjoint integer chains of step $|s|$ and is therefore a forest. Its edge-difference map on vertex values has full row rank: choose a value at one root of each component and propagate prescribed differences along the unique paths. Equivalently, its incidence matrix $B_j$ satisfies
\begin{equation}
 \rank B_j=\#\{t:j\in I_t\}.
\end{equation}
Values at distinct orbit points can be prescribed independently, and the functions $h_j$ vary independently across coordinates. Consequently
\begin{equation}\label{eq:secant-rank}
 \rank D_HE(H,x,y)=\sum_{j=1}^D\#\{t:j\in I_t\}=m.
\end{equation}
This proves the required parameter submersion without assuming that the two windows are disjoint.

\paragraph{Step 3: finite-dimensional perturbations and transversality.}
We spell out the finite-parameter argument behind Steps 1 and 2; it uses
classical parametric transversality, not a new embedding principle
\cite{thom,hirsch,golubitsky}. Let $X=SM$ for immersion or
$X=M^{(2)}$ for injectivity, and fix a base observation $H_0$.
At every $u\in X$, the constructions above give finitely many smooth
perturbations whose parameter derivatives span $\mathbb R^m$ at $u$.
Surjectivity persists in a neighborhood of $u$. Given a compact
$K\subset X$, select finitely many such neighborhoods covering $K$ and
collect their perturbations as $V_1,\ldots,V_P$. Set
\[
 H_a=H_0+\sum_{i=1}^{P}a_iV_i,\qquad a\in\mathbb R^P.
\]
The parameter derivative of $(a,u)\mapsto F(H_a,u)$ or
$(a,u)\mapsto E(H_a,u)$ is onto throughout a neighborhood of $K$.
It is independent of $a$ because the evaluation maps are linear in $H$.
Thus the total map is a submersion there. Finite-dimensional parametric
transversality gives arbitrarily small parameter vectors $a$ for which the
slice is transverse to zero. The regularity requirement is
$r>\max\{\dim X-m,0\}$. The immersion evaluation is $C^{k-1}$,
the secant evaluation is $C^k$, and both have $\dim X<m$; therefore
$k\ge2$ suffices. A transverse inverse image of zero would have dimension
$\dim X-m<0$, so it is empty. Since the perturbation family is finite,
$a\to0$ also implies $H_a\to H_0$ in $C^k$.

For immersion take $K=SM$, which is compact. For secants let
\[
 K_n=\{(x,y)\in M\times M:\operatorname{dist}(x,y)\ge1/n\}.
\]
These compact sets exhaust $M^{(2)}$. The set $\mathcal U_n$ of $H$
with $E_H\ne0$ on $K_n$ is dense by the perturbation argument. It is open
because $\min_{K_n}\|E_H\|>0$ whenever $K_n$ is nonempty, and evaluation
is continuous in the uniform norm. The space $C^k(M,\mathbb R^D)$ on
compact $M$ is Baire. Hence intersecting all $\mathcal U_n$ with the
open dense immersion set gives a dense residual set of injective
immersions. For $d=0$, omit the immersion condition and use the same
secant argument on the finite set $M^{(2)}$.

\paragraph{Step 4: inverse continuity and openness.}
An injective continuous map from compact $M$ into Euclidean space has
continuous inverse on its image: a closed subset of $M$ is compact, its
image is compact and hence closed in the image space. With immersion,
the constant-rank theorem then gives the $C^k$ inverse in image-manifold
charts. Thus the maps just obtained are $C^k$ embeddings.

We also detail why embeddings form a $C^1$-open set on compact $M$.
For an embedding $f$, at each point choose a source chart and a linear
projection of the target for which the derivative is invertible. In a
small convex chart ball, its derivative stays close enough to its value
at the center that the segment-integral formula gives a positive lower
Lipschitz bound. This persists under sufficiently small $C^1$
perturbations of $f$. Finitely many smaller balls cover $M$; a Lebesgue
number ensures all sufficiently nearby pairs lie in one such chart ball.
Pairs separated by at least that number form a compact set away from
the diagonal. Injectivity of $f$ supplies a strictly positive minimum
image separation on that set, which persists under sufficiently small
$C^0$ perturbations. Thus nearby maps remain both immersive and globally
injective. For zero-dimensional $M$, only the finite-pair separation
argument is needed. This is the usual openness result for embeddings
\cite{hirsch,golubitsky}.

For a fixed mask, $H\mapsto\Phi_{H,\mathbf I}$ is continuous from $C^k$
to $C^1$, since it uses finitely many compositions with fixed iterates of
$\varphi$. The set $\mathcal G_{\mathbf I}$ of embedding observations is
therefore open; containing a dense residual set, it is open dense.

\paragraph{Step 5: simultaneous masks and the probability clause.}
There are at most $2^{DL}$ mask words. Consequently
\[
 \mathcal G=\bigcap_{\mathbf I\in\mathcal A}\mathcal G_{\mathbf I}
\]
is a finite intersection of open dense sets and is open dense. Fix
$H\in\mathcal G$ before choosing a mask law $\nu$ supported on
$\mathcal A$. Then, explicitly,
\[
 \begin{aligned}
 \Pr_{\nu}(\text{embedding})
 &=\sum_{\mathbf I\in\mathcal A}\nu(\mathbf I)
       \mathbf 1_{\{\Phi_{H,\mathbf I}\text{ embeds}\}}\\
 &=\sum_{\mathbf I\in\mathcal A}\nu(\mathbf I)=1.
 \end{aligned}
\]
No product structure, temporal independence, or uniform sampling was used.
For a law not supported on $\mathcal A$, this argument supplies only the
lower bound $\Pr(\text{embedding})\ge\nu(\mathcal A)$; a word outside
$\mathcal A$ might still embed. If $\mathcal A$ is empty, the simultaneous
assertion is vacuous and there is no probability law supported on it.
\end{proof}

\paragraph{Why the hypotheses cannot simply be omitted.}
For $M=S^1$, $\varphi=\mathrm{id}$, $D=2$, $p=1$, and $L=3$,
observe only the first sensor at each epoch. Then $Lp=3=2d+1$, but
$\Phi(x)=(h_1(x),h_1(x),h_1(x))$. No continuous map $h_1:S^1\to\mathbb R$
is injective, so the count condition alone cannot establish the theorem.
Likewise, a constant prescribed observation map fails even for an
aperiodic rotation: genericity is not a guarantee for every sensor.
If the all-erased word has positive probability and $M$ has two distinct
points, embedding fails on that event. These examples explain the roles
of the dynamics, generic-observation, and supported-mask hypotheses;
they do not establish that the sufficient exclusion range or scalar
count is necessary or sharp. The theorem also does not identify a
learned numerical estimator with its ideal inverse.

\subsection{Proof of Corollary~\ref{cor:burst}}\label{app:proof-cor-burst}
\begin{proof}
For (i), bounded-burst words form a subset of the fixed-budget words already covered by Theorem~\ref{thm:main}. For (ii), partition the first $(K+1)\lfloor L/(K+1)\rfloor$ epochs into disjoint blocks of length $K+1$. Each coordinate is observed at least once in every block; otherwise that block would contain $K+1$ consecutive erasures of that coordinate. Summing over blocks and coordinates gives \eqref{eq:burst-count}. The variable-budget conclusion of Theorem~\ref{thm:main} applies to every word satisfying \eqref{eq:burst-bound}, simultaneously. Explicitly, with $m=\lfloor L/(K+1)\rfloor$,
\[
 \sum_{t=0}^{L-1}|I_t|
 =\sum_{j=1}^D\sum_{t=0}^{L-1}{\bf1}_{\{j\in I_t\}}
 \ge\sum_{j=1}^D\sum_{b=0}^{m-1}1=Dm.
\]
The unpartitioned final epochs contribute a nonnegative count. This also
covers $m=0$; no positive lower count can then be inferred from the burst
constraint alone.
\end{proof}

\subsection{Proof of Proposition~\ref{prop:periodic}}\label{app:proof-prop-periodic}
\begin{proof}
Choose an eigenbasis $v_1,\ldots,v_d$ of the complexified tangent space, with $Av_k=\lambda_kv_k$. Put
\[
 \alpha_k=Dh_{j_*}(\varphi^{r_0}x_0)D(\varphi^{r_0})_{x_0}v_k.
\]
The selected derivative rows in this basis form a square matrix with entries $\alpha_k\lambda_k^{q_0+i}$. Its determinant is
\begin{equation}\label{eq:vandermonde}
 \left(\prod_{k=1}^d\alpha_k\lambda_k^{q_0}\right)
 \prod_{1\le a<b\le d}(\lambda_b-\lambda_a).
\end{equation}
All eigenvalues are nonzero and distinct. Each condition $\alpha_k=0$ is a proper real-linear subspace of the real covector space (for a nonreal eigenvector it represents the corresponding real conditions). A finite union of such subspaces has empty interior. Thus $\alpha_k\ne0$ for all $k$ is open dense and is realizable by local perturbation of $h_{j_*}$. The determinant is then nonzero. Complex full rank of a real derivative matrix is equivalent to real full rank.
\end{proof}

\subsection{Proof of Proposition~\ref{prop:probability}}\label{app:proof-prop-probability}
\begin{proof}
For the chosen $H$, every word with $N_L\ge2d+1$ is an embedding word. Thus failure implies $N_L\le2d$. To show the concentration step explicitly, put $n=LD$ and $a=\mu-2d>0$.
Hoeffding's bounded-variable exponential estimate gives, for $\lambda>0$,
\[
 \mathbb E e^{-\lambda(N_L-\mu)}\le e^{n\lambda^2/8}.
\]
Markov's inequality therefore gives
\[
 \mathbb P(N_L-\mu\le-a)
 \le e^{-\lambda a+n\lambda^2/8}.
\]
Choosing $\lambda=4a/n$ yields $e^{-2a^2/n}$, which is
\eqref{eq:hoeffding} \cite{hoeffding}. The first inequality is a containment, not an equality: words with fewer observations may still embed.
\end{proof}

\section{Failure of a measure-zero periodic-set hypothesis}
Take $M=\mathbb T^2$, $D=3$, $p=1$, and $L=5$. For $0<\varepsilon<1$, set
\[
 \begin{gathered}
 f(\psi)=\psi+\varepsilon(1-\cos\psi)\pmod{2\pi},\\
 \varphi(\theta,\psi)=(\theta,f(\psi)).
 \end{gathered}
\]
Since $f'(\psi)=1+\varepsilon\sin\psi>0$, this is a smooth diffeomorphism. On the lift interval $(0,2\pi)$, $\psi<f(\psi)<2\pi$, so there are no periodic points there. Thus $\Fix(\varphi^q)=\mathbb T^1\times\{0\}$ for every $q\ge1$, a set of zero two-dimensional volume.

If only coordinate 1 is observed at every epoch, then on this fixed circle
\[
 \Phi_{H,\mathbf I}(\theta,0)=(h_1(\theta,0),\ldots,h_1(\theta,0)).
\]
This is not injective, since a circle cannot be continuously embedded in a line. Also $D\varphi_{(\theta,0)}=I_2$, so the five derivative rows coincide and the rank is at most one. This failure holds for every $H$, although $Lp=2d+1$. The short-period exclusion removes this obstruction; it is a substantive restriction of the theorem's domain.

\section{Proofs of the forecast error bounds}\label{app:forecast-proofs}
This appendix gives the complete proofs of the ten forecast and inverse-error statements in Section~\ref{sec:forecast} and Appendix~\ref{app:forecast-cases}. All formal results use one automatic numbering sequence; the mask-transition proof appears in Appendix~\ref{app:new:prop:mask-transition}.

\subsection{Proof of Lemma~\ref{fc:lem:inverse}}\label{app:proof-fc-lem-inverse}
\begin{proof}
Suppose distinct $x_j,y_j\in M$ have secant ratios tending to zero. Compactness gives a subsequence with $x_j\to x$ and $y_j\to y$. If $x\ne y$, continuity yields $\Psi(x)=\Psi(y)$, contradicting injectivity. Thus $x=y=p$. A $C^1$ graph chart at $p$ and the fundamental theorem of calculus in the chart imply, after taking a further subsequence,
\[
 \frac{x_j-y_j}{\norm{x_j-y_j}}\longrightarrow v\in T_pM,\qquad \norm v=1,
\]
and
\[
 \frac{\Psi(x_j)-\Psi(y_j)}{\norm{x_j-y_j}}\longrightarrow D\Psi_p v.
\]
For completeness, in graph coordinates $q(u)=(u,\gamma(u))$, both differences equal their derivatives at the limiting chart point applied to $u_j-v_j$, plus $o(\norm{u_j-v_j})$. The graph parametrization is locally bi-Lipschitz, so these remainders are also $o(\norm{x_j-y_j})$. This proves the two displayed limits. Immersion gives $D\Psi_pv\ne0$, again a contradiction. The inverse inequality and the finite-minimum assertion follow immediately.
\end{proof}

\subsection{Proof of Corollary~\ref{fc:cor:recovery}}\label{app:proof-fc-cor-recovery}
\begin{proof}
Insert $\widetilde\Psi(\widehat x_0)$ and $y$ between $\Psi(\widehat x_0)$ and $\Psi(x_0)$. The triangle inequality bounds their difference by $\beta+r_{\rm fit}+\eta$; apply Lemma~\ref{fc:lem:inverse}. In the least-squares case, the feasible true state has residual at most $\eta$, so a global minimizer has no larger residual.
\end{proof}

\subsection{Proof of Theorem~\ref{fc:thm:continuous}}\label{app:proof-fc-thm-continuous}
\begin{proof}
Set $e=\widehat x-x$ and
\[
 d(t)=\widehat f(t,\widehat x(t))-f(t,\widehat x(t))+r_{\rm num}(t).
\]
Then $\dot e=f(t,\widehat x)-f(t,x)+d$ and $\norm d\le b$. Hence
\begin{equation}
 \tfrac12\frac{d}{dt}\norm e^2\le\ell(t)\norm e^2+b(t)\norm e
 \quad\text{almost everywhere}.\label{fc:eq:energy}
\end{equation}
To handle zeros of $e$ and negative $\ell$ rigorously, fix $a>0$ and set $v_a=(\norm e^2+a^2)^{1/2}$. Since $v_a\ge a$,
\begin{align*}
 \dot v_a&\le\ell(t)\frac{\norm e^2}{v_a}+b(t)\frac{\norm e}{v_a}\\
 &\le\ell(t)v_a+b(t)+|\ell(t)|a.
\end{align*}
Multiply by $e^{-A(t)}$ and integrate:
\[
 v_a(t)\le e^{A(t)}v_a(0)+
 \int_0^t e^{A(t)-A(s)}\big(b(s)+a|\ell(s)|\big)\,ds.
\]
Integrability on the finite interval permits $a\downarrow0$, proving \eqref{fc:eq:general}.
\end{proof}

\subsection{Proof of Corollary~\ref{fc:cor:pipeline}}\label{app:proof-fc-cor-pipeline}
\begin{proof}
Substitute \eqref{fc:eq:recovery} in \eqref{fc:eq:general}. For the constant-bound claim, use $A(t)-A(s)\le\ell(t-s)$ and integrate $e^{\ell(t-s)}$.
\end{proof}

\subsection{Proof of Corollary~\ref{fc:cor:horizon}}\label{app:proof-fc-cor-horizon}
\begin{proof}
For $\ell>0$, rewrite $B(t)=e^{\ell t}(\delta_0+b/\ell)-b/\ell$. This is increasing unless both error bounds vanish. Solving $B(t)=\Delta$ proves \eqref{fc:eq:horizon}. For $\ell=0$ the expression is linear. For $\ell=-a$, $B(t)$ is a convex combination of $\delta_0$ and $s$. If $s>\Delta\ge\delta_0$, it increases to $s$ and solving its threshold equation gives \eqref{fc:eq:contracthorizon}.
\end{proof}

\subsection{Proof of Proposition~\ref{fc:prop:ag}}\label{app:proof-fc-prop-ag}
\begin{proof}
Fix $t$ and differentiate $q(s)=\Phi^{t-s}(\widehat x(s))$ almost everywhere. The flow identity implies $D\Phi^u(z)f(z)=f(\Phi^u(z))$. Therefore
\[
 \begin{aligned}
 q'(s)&=-f(\Phi^{t-s}(\widehat x(s)))\\
 &\quad+D\Phi^{t-s}(\widehat x(s))\dot{\widehat x}(s)\\
 &=D\Phi^{t-s}(\widehat x(s))d(s).
 \end{aligned}
\]
Integrate from $0$ to $t$, using $q(t)=\widehat x(t)$ and $q(0)=\Phi^t(\widehat x_0)$. Subtract $\Phi^t(x_0)$. The initial term is bounded by integrating $D\Phi^t$ along the initial segment; taking norms gives \eqref{fc:eq:kernel}.
\end{proof}

\subsection{Proof of Theorem~\ref{fc:thm:discrete}}\label{app:proof-fc-thm-discrete}
\begin{proof}
Insert $F_n(\widehat x_n)$ and use the triangle inequality to obtain $E_{n+1}\le q_nE_n+b_n$. Repeated substitution proves \eqref{fc:eq:discrete} by induction. Summing the geometric series proves \eqref{fc:eq:geometric}.
\end{proof}

\subsection{Proof of Theorem~\ref{fc:thm:collocation}}\label{app:proof-fc-thm-collocation}
\begin{proof}
The true step is $x_{n+1}=x_n+\tfrac12[g(x_n)+g(x_{n+1})]+d_g(x_n)$. Subtract this from \eqref{fc:eq:forecastcoll}, take norms, and use the Lipschitz condition:
\[
 E_{n+1}\le E_n+\tfrac\kappa2(E_n+E_{n+1})+\rho_n+\eps_{\rm col}.
\]
Since $1-\kappa/2>0$, rearrangement proves \eqref{fc:eq:colrec}. Apply Theorem~\ref{fc:thm:discrete}, noting that $q-1=\kappa/(1-\kappa/2)$, to obtain \eqref{fc:eq:colbound}.
\end{proof}

\subsection{Proof of Proposition~\ref{fc:prop:coverage}}\label{app:proof-fc-prop-coverage}
\begin{proof}
For a true training pair, subtract its defect from the reconstructed residual. The difference has norm at most
\[
 (1+\kappa/2)\big(\norm{\widetilde z_i-x_i}+
 \norm{\widetilde z_i^+-F(x_i)}\big)\le(2+\kappa)\zeta.
\]
Thus $\norm{d_g(x_i)}\le R_{\max}+(2+\kappa)\zeta$. Also, for $x,y\in K$,
\begin{align*}
 \norm{d_g(x)-d_g(y)}
 &\le\norm{F(x)-F(y)}+\norm{x-y}\\
 &\quad+\tfrac\kappa2\big(\norm{x-y}+\norm{F(x)-F(y)}\big)\\
 &\le(1+q_F)(1+\kappa/2)\norm{x-y}.
\end{align*}
Choose a training point within $r_{\rm cov}$ of $x$ and apply this estimate.
\end{proof}

\subsection{Proof of Corollary~\ref{fc:cor:trapezoid}}\label{app:proof-fc-cor-trapezoid}
\begin{proof}
For $v(s)=x'(s)$ on $[0,h]$, integration by parts twice yields
\[
 \int_0^h v(s)\,ds-\tfrac h2\big(v(0)+v(h)\big)
 =-\tfrac12\int_0^h s(h-s)v''(s)\,ds.
\]
Its norm is at most $M_3h^3/12$. Since $v(s)=f(x(s))$, this is the defect of the true step using the increment field $hf$. Replacing $hf$ at its two endpoints by $g$ adds at most $h\eps_f$. This proves \eqref{fc:eq:quaddefect}. Apply Theorem~\ref{fc:thm:collocation} with the valid Lipschitz upper bound $\kappa=hL_g$.
\end{proof}

\section{Proofs of the mask-dependent extensions}\label{app:mask-extensions}

\subsection{Proof of result~\ref{thm:periodic}: Exact generic local-immersion criterion}\label{app:new:thm:periodic}
\begin{proof}
Write $x_r=\varphi^r(x_0)$ and $P_r=D(\varphi^r)_{x_0}$. The chain rule gives
\[
 D(h_j\circ\varphi^{r+qp_0})_{x_0}=Dh_j(x_r)P_rA^q.
\]
The covectors $c_g=Dh_j(x_r)P_r$ vary independently: different phases are distinct points, different coordinates are independent scalar functions, and every $P_r$ is invertible. Bump functions in disjoint neighborhoods prescribe these finite first jets and give a continuous finite-dimensional right inverse to their evaluation map.

Choose an eigenbasis $v_k$ of the complexified tangent space and put $\alpha_{gk}=c_gv_k$. The derivative in this basis is
\[
 \mathcal O(\alpha)=\begin{bmatrix}
 V_{g_1}\diag(\alpha_{g_1})\\[-1mm]
 V_{g_2}\diag(\alpha_{g_2})\\
 \vdots
 \end{bmatrix}.
\]
Choose $d$ rows, with row sets $R_g$ in the groups. Expanding their determinant by block rows yields
\begin{equation}\label{eq:minor}
 \sum_{(S_g)}\operatorname{sgn}(S_g)
 \prod_g\left[\det V_g[R_g,S_g]\prod_{k\in S_g}\alpha_{gk}\right],
\end{equation}
where the sum runs over column partitions with $|S_g|=|R_g|$. Different partitions give different monomials in the variables $\alpha_{gk}$, so their coefficients cannot cancel. Some $d$-row determinant is a nonzero polynomial if and only if a partition as in item 3 exists.

This argument remains valid for real covectors when eigenvalues occur in conjugate pairs. Complexifying the real covector coordinates gives an invertible complex-linear change of variables to $\alpha$. A nonzero polynomial after this change cannot vanish on the whole real covector space: a polynomial vanishing on $\R^n$ is identically zero by successive one-variable restriction. Since the derivative is real in a real tangent basis, complex full rank is equivalent to real full rank. Thus at least one real minor is a nonzero real polynomial whenever item 3 holds.

Its nonvanishing set is open dense in finite real jet space. The bump-function right inverse transfers density to the observation-function space, and rank continuity gives openness. This proves $3\Rightarrow2\Rightarrow1$. If item 3 fails, every $d$-row determinant is identically zero by \eqref{eq:minor}, proving $1\Rightarrow3$ and the final assertion.

For $3\Leftrightarrow4$, each $V_g$ defines its column linear matroid on the common eigenmode set. The equivalence is the matroid-partition theorem of Edmonds and Fulkerson \cite[Theorem 1c, p.~150]{EF}. This combinatorial theorem is explicitly credited; the determinant proof of items 1--3 does not require it.
\end{proof}

\subsection{Proof of result~\ref{new:periodic:2}: Consecutive-cycle sufficient conditions}\label{app:new:new:periodic:2}
\begin{proof}
For $d$ consecutive exponents the determinant is
\[
 \left(\prod_{k=1}^d\lambda_k^{q_0}\right)
 \prod_{a<b}(\lambda_b-\lambda_a)\ne0.
\]
For a consecutive block of length $n_g$, any subset of at most $n_g$ columns is independent by the same Vandermonde argument, using the first corresponding number of rows. Thus one can distribute the $d$ modes among groups if their total capacity is at least $d$. If total capacity is smaller than $d$, the entire observation matrix has fewer than $d$ rows.
\end{proof}

\subsection{Proof of result~\ref{cor:fixed-periodic}: Arbitrary switching with a fixed per-epoch budget}\label{app:new:cor:fixed-periodic}
\begin{proof}
Necessity follows because the derivative has only $bL$ rows. For sufficiency, fix an eigenmode subset $S$ of cardinality $s$. Phase $r$ occurs $n_r$ times, at cycle indices $0,\ldots,n_r-1$. The full $n_r$-by-$s$ Vandermonde matrix has rank $\min(n_r,s)$, so choose a set $B_r$ of that many independent rows. Every cycle index occurs in exactly $b$ coordinate groups. For each coordinate $j$, the rows $B_r\cap Q_{(j,r)}$ are independent, giving
\[
 \sum_j r_{(j,r)}(S)\ge\sum_j|B_r\cap Q_{(j,r)}|=b\min(n_r,s).
\]
If some $n_r\ge s$, summing over phases is at least $bs\ge s$. Otherwise it equals at least $b\sum_r n_r=bL\ge d\ge s$. Thus every inequality in \eqref{eq:matroid} holds, and Theorem~\ref{thm:periodic} applies. A finite intersection yields one open dense good set for all fixed-budget masks.
\end{proof}

\subsection{Proof of result~\ref{cor:local-reliability}: Exact finite-window local reliability}\label{app:new:cor:local-reliability}
\begin{proof}
Intersect the finitely many open dense success sets for masks in $\mathcal G$. Masks outside $\mathcal G$ fail for every $H$. If $\mathcal G$ is empty, use the entire observation space as $\mathcal U$.
Otherwise the finite intersection is open dense. For every fixed $H\in\mathcal U$,
\[
 \begin{split}
 \Prob_\pi\{\rank D\Phi_{H,\mathcal I}(x_0)=d\}
 &=\sum_{I\in\mathcal F}\pi(I)
 {\bf1}_{\{\rank D\Phi_{H,I}(x_0)=d\}}\\
 &=\sum_{I\in\mathcal G}\pi(I).
 \end{split}
\]
The good observation set is chosen before $\pi$; this is not a probability
law on observation functions. This proves \eqref{eq:localprob}.
\end{proof}

\subsection{Proof of result~\ref{new:horizon:5}: Finite horizons versus all time}\label{app:new:new:horizon:5}
\begin{proof}
The fixed-count assertion follows pathwise from the assumed common good observation law. The finite-horizon assertion is the union bound. In the last model the probability of an empty block is
\[
 q_0=\prod_{j=1}^D(1-\rho_j)^L>0.
\]
Disjoint blocks are independent. The probability that none of the next $n$ blocks is empty equals $(1-q_0)^n\to0$; applying this after each finite block index and taking a countable intersection shows infinitely many empty blocks almost surely. An empty observation map is constant and cannot embed a space with at least two points.
\end{proof}

\subsection{Proof of result~\ref{thm:geometry-mask}: Geometry-aware random-mask lower bound}\label{app:new:thm:geometry-mask}
\begin{proof}
Fix $v\in\mathcal S$. Put $a_{t,j}=v_{t,j}^2$ and $b=\mu/N$, so $0\le a_{t,j}\le b$ and $\sum_{t,j}a_{t,j}=1$. Write
\[
 Z=\|P_{\mathbf I}v\|^2=\sum_{t,j}s_{t,j}a_{t,j}.
\]
We first establish, for every $\lambda\ge0$,
\begin{equation}\label{eq:mask-mgf}
 \mathbb E e^{-\lambda Z}
 \le \exp\left(\frac{q}{b}(e^{-\lambda b}-1)\right).
\end{equation}
For independent Bernoulli coordinates,
\[
\begin{aligned}
 \mathbb E e^{-\lambda Z}
 &=\prod_{t,j}\{1-q+qe^{-\lambda a_{t,j}}\}\\
 &\le\exp\left(q\sum_{t,j}(e^{-\lambda a_{t,j}}-1)\right).
\end{aligned}
\]
For exact-$p$ sampling, indicators within an epoch are dependent and must not be treated as independent. Set $u_j=e^{-\lambda a_{t,j}}$. The elementary symmetric-mean inequality gives
\[
\begin{aligned}
 \mathbb E\exp\left(-\lambda\sum_{j\in I_t}a_{t,j}\right)
 &=\frac{e_p(u_1,\ldots,u_D)}{\binom Dp}\\
 &\le\left(\frac1D\sum_{j=1}^D u_j\right)^p\\
 &\le\exp\left(q\sum_{j=1}^D(u_j-1)\right).
\end{aligned}
\]
For completeness, the symmetric-mean inequality follows by maximizing $e_p$ at a fixed sum of nonnegative $u_j$. For a pair $u_i,u_j$, the polynomial has the form $A u_i u_j+B(u_i+u_j)+C$, where $A\ge0$; replacing the pair by its average does not decrease it. Among maximizers on the compact simplex, choose one minimizing $\sum_j u_j^2$. If two entries differ, averaging them does not decrease $e_p$ and strictly decreases this sum of squares, a contradiction. Thus the maximizing vector has equal coordinates, proving the bound. The case $p=1$ is equality. Independence between epochs permits multiplication of the block bounds.

In either law, convexity on $[0,b]$ yields
\[
 e^{-\lambda a}-1\le \frac{a}{b}(e^{-\lambda b}-1).
\]
Sum this inequality and use $\sum a_{t,j}=1$ to obtain \eqref{eq:mask-mgf}. For $0<\delta<1$, Markov's inequality gives
\[
 \mathbb P\{Z\le q(1-\delta)\}
 \le\inf_{\lambda\ge0}
 \exp\left(\lambda q(1-\delta)+\frac q b(e^{-\lambda b}-1)\right).
\]
Choose $\lambda=-b^{-1}\log(1-\delta)$. The exponent becomes
\[
 -\frac q b\{\delta+(1-\delta)\log(1-\delta)\}
 \le-\frac{q\delta^2}{2b}.
\]
The last inequality follows because the derivative of the expression in braces is $-\log(1-\delta)\ge\delta$ and its value at zero is zero. Thus
\begin{equation}\label{eq:fixed-secant-tail}
 \mathbb P\{\|P_{\mathbf I}v\|^2\le q(1-\delta)\}
 \le\exp\left(-\frac{qN\delta^2}{2\mu}\right).
\end{equation}

Choose a deterministic $\eta$-net $\mathcal V\subset\mathcal S$, where $\eta=\epsilon\sqrt q/4$, and apply \eqref{eq:fixed-secant-tail} with $\delta=\epsilon/2$ to all its points. Except on an event of probability at most
\[
 |\mathcal V|\exp\left(-\frac{qN\epsilon^2}{8\mu}\right),
\]
every $u\in\mathcal V$ satisfies $\|P_{\mathbf I}u\|\ge\sqrt{q(1-\epsilon/2)}$. For any $v\in\mathcal S$, choose $u\in\mathcal V$ with $\|u-v\|\le\eta$. Since $\|P_{\mathbf I}\|\le1$,
\[
 \|P_{\mathbf I}v\|\ge\sqrt q\left(\sqrt{1-\epsilon/2}-\epsilon/4\right)
 \ge\sqrt{q(1-\epsilon)}.
\]
Indeed,
\[
 \sqrt{1-\epsilon/2}-\sqrt{1-\epsilon}
 =\frac{\epsilon/2}{\sqrt{1-\epsilon/2}+\sqrt{1-\epsilon}}
 \ge\epsilon/4.
\]
Apply this to each normalized secant, multiply by $\|F(x)-F(y)\|$, and use \eqref{eq:full-secant-lower}. This proves the claimed event and bound.
\end{proof}

\subsection{Proof of result~\ref{cor:nonuniform-geometry}: Nonuniform survival and variable exact budgets}\label{app:new:cor:nonuniform-geometry}
\begin{proof}
The block or Bernoulli calculation yields
\[
\begin{aligned}
 \mathbb E e^{-\lambda Z}\\
 &\le\exp\left(\sum_{t,j}\rho_{t,j}(e^{-\lambda a_{t,j}}-1)\right)\\
 &\le\exp\left(\frac{e^{-\lambda b}-1}{b}
                 \sum_{t,j}\rho_{t,j}a_{t,j}\right).
\end{aligned}
\]
Since $e^{-\lambda b}-1\le0$ and the weighted sum is at least $\vartheta$, the last expression is at most $\exp[\vartheta(e^{-\lambda b}-1)/b]$. Repeat the proof of the theorem.
\end{proof}

\subsection{Proof of result~\ref{cor:noise-rolling}: Noise amplification and finitely many rolling windows}\label{app:new:cor:noise-rolling}
\begin{proof}
The triangle inequality bounds $\|P_{\mathbf I}(F(\widehat x)-F(x))\|$ by $r+\beta+\eta_0$. Divide by the lower constant in \eqref{eq:mask-stable-lower}. The exact-map least-squares comparison uses the feasible true state. For rolling windows, let $A_w$ be the good geometry event in window $w$.
If the common bound $B_\epsilon$ applies to each window, then
\[
 \Prob\Bigl(\bigcap_{w=1}^W A_w\Bigr)
 =1-\Prob\Bigl(\bigcup_{w=1}^W A_w^c\Bigr)
 \ge1-\min\{1,WB_\epsilon\}.
\]
For different certified geometries or marginals, replace $WB_\epsilon$
by $\sum_w B_{\epsilon,w}$. No factorization of probabilities is used.
The noise, residual and map-approximation hypotheses must separately
hold in each window; the union bound only certifies the geometry events.
\end{proof}

\subsection{Proof of result~\ref{prop:mask-transition}: Mask-aware transition identity and stable forcing}\label{app:new:prop:mask-transition}
\begin{proof}
Successive inverse pairs $\Psi_{I_k}^{-1}\circ\Psi_{I_k}$ cancel on $M$, proving \eqref{eq:switch-telescope}; composition of the endpoint Lipschitz bounds proves \eqref{eq:switch-lipschitz}. For the second assertion,
\[
\begin{aligned}
 \|\widehat x_{k+1}-\varphi(\widehat x_k)\|\\
 &\le\frac{\|\Psi_{I_{k+1}}(\widehat x_{k+1})-
                   \Psi_{I_{k+1}}(\varphi(\widehat x_k))\|}{c_{I_{k+1}}}\\
 &\le\frac{r_k}{c_{I_{k+1}}}.
\end{aligned}
\]
Thus $E_{k+1}\le Q E_k+r_k/c_{I_{k+1}}$. Induction yields \eqref{eq:switch-forcing}. The mask sequence may be arbitrary, including a realized random sequence; this algebra requires no mask independence.
\end{proof}
\clearpage

\section{Quantitative certificates and classical sampling specialization}\label{app:quantitative-certificates}
\subsection{Proof of Corollary~\ref{cor:subspace-certificate}}\label{app:subspace-certificate}
\begin{proof}
Write $a_{t,j}=Q^\top e_{t,j}$ and
$A_{t,j}=a_{t,j}a_{t,j}^{\top}$. The sampled Gram matrix is
$G=\sum_t\sum_{j\in I_t}A_{t,j}$, each $A_{t,j}$ is positive semidefinite
with $\lambda_{\max}(A_{t,j})\le b$, and $\mathbb EG=qI_r$.

For a fixed Hermitian matrix $C$ and $\theta>0$, the function
$X\mapsto\operatorname{tr}\exp[-\theta(C+X)]$ is convex.
The matrix-valued convex comparison for sampling without replacement
\cite[Sec.~D]{GN} therefore bounds this expectation for one
exact-budget epoch by the expectation for $p$ independent uniform draws
with replacement from its indexed population $\{A_{t,j}:1\le j\le D\}$.
Condition on the other epochs and apply this comparison successively to
each independent epoch. Thus
\[
 \mathbb E\operatorname{tr}e^{-\theta G}
 \le \mathbb E\operatorname{tr}e^{-\theta\widetilde G},
\]
where $\widetilde G$ is the sum of $Lp$ independent positive semidefinite
draws with eigenvalues bounded by $b$ and total expectation $qI_r$.
For the lower tail, $\lambda_{\min}(G)\le q(1-\epsilon)$ implies $\operatorname{tr}e^{-\theta G}\ge e^{-\theta q(1-\epsilon)}$. Markov's inequality and the trace comparison reduce the tail to the independent model. The matrix Laplace argument underlying
\cite[Cor.~5.2 and Rem.~5.3]{Tropp2012} now gives
$r\exp[-(q/b)h_-(\epsilon)]$ for the stated lower-tail event. Explicitly,
the independent trace bound is $r\exp[(q/b)(e^{-\theta b}-1)]$.
After Markov's inequality, the exponent is
$\theta q(1-\epsilon)+(q/b)(e^{-\theta b}-1)$.
Substitution of $\theta=-b^{-1}\log(1-\epsilon)$ gives
$-(q/b)h_-(\epsilon)$.
Since $h_-(\epsilon)\ge\epsilon^2/2$, the simpler bound follows.
Bernoulli sampling permits direct use of matrix Chernoff.
For a vector $v=Qz\in U$, the good event implies
$\|P_{\mathbf I}v\|^2=z^\top Gz\ge q(1-\epsilon)\|v\|^2$.
Apply this to every secant and use the assumed full-map lower bound.
\end{proof}
\subsection{Circle geometry, conditioning, and noise}\label{app:circle-certificate}
Let $M=S^1\subset\mathbb R^2$, $\varphi(x)=R_\omega x$, and $H(x)=x$, where
$R_\omega$ is planar rotation through $\omega=2\pi\alpha$ and
$\alpha=(\sqrt5-1)/2$. This is a smooth, invertible dynamics on a compact
one-dimensional state space. Its irrational rotation has no periodic points.
For a terminal window of length $L$, write
\[
 A_L=\begin{pmatrix}R_{-(L-1)\omega}\\\vdots\\R_{-\omega}\\I_2\end{pmatrix},
 \qquad F_L(x)=A_Lx\in\mathbb R^{2L}.
\]
At each epoch retain exactly one of the two physical coordinates, each with
probability $1/2$, independently across epochs. Thus $D=2$, $p=1$, $q=1/2$,
$N=2L$, and the realized data consist of precisely $L$ scalar observations.
The mask, including the identity of each observed coordinate, is known.
No replicated sensor channels are introduced. The different delay blocks
rotate the measurement directions; irrationality ensures that directions
from distinct blocks do not coincide even modulo sign.

\paragraph{Certified circle geometry.}
For every positive integer $L$, the full lower-secant constant is exactly
$c_F=\sqrt L$, the coordinate coherence is exactly $\mu=2$, and, for
$0<\eta<2$, the normalized-secant covering number, using closed balls (distance $\le\eta$), is
\[
 \mathcal N(\mathcal S,\eta)
 =\left\lceil\frac{\pi}{2\arcsin(\eta/2)}\right\rceil.
\]
In particular, with
\[
 m_\epsilon=\left\lceil
 \frac{\pi}{2\arcsin(\epsilon/(8\sqrt2))}\right\rceil,
\]
\[
 \overline B_\epsilon=\min\{1,m_\epsilon e^{-L\epsilon^2/16}\},
\]
with probability at least $1-\overline B_\epsilon$, every pair $x,y\in S^1$
satisfies
\[
 \|P_{\mathbf I}F_L(x)-P_{\mathbf I}F_L(y)\|
 \ge \sqrt{L(1-\epsilon)/2}\,\|x-y\|.
\]

\begin{proof}
Orthogonality gives $A_L^\top A_L=L I_2$, so all full secant distances are
multiplied by $\sqrt L$. Every direction $u\in S^1$ occurs as a state secant,
by choosing the antipodal pair $x=u,y=-u$. Consequently
$\mathcal S=\{A_Lu/\sqrt L:u\in S^1\}$, an isometric copy of the unit circle.
Every row of $A_L$ has norm one. Therefore
$\sup_{u\in S^1}\|A_Lu/\sqrt L\|_\infty^2=1/L$ and $\mu=2L/L=2$.

A Euclidean ball of radius $\eta$ centered on the unit circle covers a
closed arc of angular length $4\arcsin(\eta/2)$. Covering its full angular
length $2\pi$ therefore requires at least
$\lceil\pi/[2\arcsin(\eta/2)]\rceil$ centers. Equally spaced centers attain
this number: the largest angle to a nearest center is $\pi/m$, giving
maximum chord distance $2\sin(\pi/(2m))\le\eta$. Isometry transfers this
cover to $\mathcal S$. Substitution of $q=1/2$, $N=2L$, $\mu=2$, and
$\eta=\epsilon/(4\sqrt2)$ into Theorem~\ref{thm:geometry-mask} proves the
claim. All constants and the net are fixed before drawing a mask.
\end{proof}

For example, $L=512$ and $\epsilon=1/2$ give $m_\epsilon=36$ and
\[
 \overline B_{1/2}=36e^{-8}<0.01208,
 \qquad c_{\mathbf I}\ge\sqrt{128}.
\]
Thus the all-state lower bound holds with probability greater than
$0.98792$ using two physical sensors and one retained scalar per epoch.
Shorter, weaker certificates are also possible: $(L,\epsilon)=(128,0.9)$
gives $20e^{-6.48}<0.03068$ and $c_{\mathbf I}\ge\sqrt{6.4}$;
$(192,0.75)$ gives $24e^{-6.75}<0.02811$ and
$c_{\mathbf I}\ge\sqrt{24}$. These are per-window guarantees, not a
uniform guarantee over an indefinitely long mask sequence. For $W$ windows,
a valid union bound is $\min\{1,W\overline B_\epsilon\}$, which may again
be vacuous. Acquiring a longer window also has an explicit cost in data
collection and assumes a fixed, known dynamics over that window.

\paragraph{A shorter-window comparison using classical matrix Chernoff.}
The general Theorem~\ref{thm:geometry-mask} secant-net argument is deliberately not optimized for this
known two-dimensional linear subspace. For this example, classical matrix
Chernoff bounds give a sharper alternative, without changing Theorem~\ref{thm:geometry-mask} or claiming
a new principle. Let $a_t^\top$ be the retained row at epoch $t$ and set
$X_t=a_ta_t^\top$. These independent positive semidefinite $2\times2$
matrices satisfy $\lambda_{\max}(X_t)=1$ and
$\mathbb E X_t=I_2/2$. Thus Tropp's matrix Chernoff bound
\cite[Corollary~5.2 and Remark~5.3]{Tropp2012} gives
\[
\begin{aligned}
 &\Pr\{c_{\mathbf I}^2\le L(1-\epsilon)/2\}\\
 &\quad\le \min\{1,\,2e^{-(L/2)h(\epsilon)}\}
 \\
 &\quad\le \min\{1,\,2e^{-L\epsilon^2/4}\}.
\end{aligned}\]
\[h(\epsilon)=\epsilon+(1-\epsilon)\log(1-\epsilon).
\]
The identity $c_{\mathbf I}^2=\lambda_{\min}(\sum_t X_t)$ follows from
the complete circle secant set established above; the last inequality uses
$h(\epsilon)\ge\epsilon^2/2$. In particular, at $L=64$ and
$\epsilon=1/2$, the squared lower threshold is $16$ and the simpler failure
bound is $2e^{-4}<0.03664$. The entropy version is
\[
 2\exp[-32(1/2+(1/2)\log(1/2))]<0.014751.
\]
By contrast, the Theorem~\ref{thm:geometry-mask} expression at these same parameters is
$\min\{1,36e^{-1}\}=1$. This comparison reports a useful refinement
available from an existing matrix concentration theorem; it does not
retroactively strengthen the Theorem~\ref{thm:geometry-mask} formula. The $L=512$ example
above remains a direct nonvacuity demonstration of that original formula.
For $L=64$ the total-window-noise gain is $2/4=1/2$, whereas the
coordinate-bounded-noise gain remains $4$, under the respective noise
conventions below. Since $p=1$ here, the retained row matrices are directly
independent across epochs; no independence of coordinate indicators within
an epoch is assumed.

\paragraph{Geometric and physical meaning.}
The full delay map distributes every state-space direction evenly in total
energy: its two singular values are both $\sqrt L$. Coherence $\mu=2$
certifies that no single scalar measurement can carry more than $1/L$ of
any normalized secant's energy. Independent per-epoch coordinate choices
then preserve a uniformly positive fraction of every secant's energy,
with the displayed high probability. The lower bound rules out flattening
of the observed circle in any direction, not just collisions among finitely
many sampled states. The length $L$ supplies measurements of the same
terminal state through a known invertible dynamics, not $L$ independent
unknown states. Irrationality is used to make the example dynamically
nonperiodic and its directions distinct; the proof of the Theorem~\ref{thm:geometry-mask} constants
itself uses only orthogonality.

\paragraph{Why the guarantee is probabilistic despite all-mask injectivity.}
For $L\ge2$, every mask in this example is injective: two retained rows
from distinct epochs could be parallel only if
$(t-s)\omega$ were an integer multiple of $\pi/2$, contrary to the
irrationality of $\omega/(2\pi)$. Their rank is therefore two, which
also certifies immersion on the circle. Nevertheless, the particular
$\epsilon=1/2$ lower bound $c_{\mathbf I}^2\ge L/4$ cannot hold for all
masks. To see this without numerical evidence, choose a unit direction
$u=(\cos\theta,\sin\theta)$ and, at every epoch, retain the coordinate
whose squared response in that direction is smaller. Averaging the total
retained energy over uniform $\theta$ gives
\[
 L\,\frac{1}{2\pi}\int_0^{2\pi}
 \min\{\cos^2\theta,\sin^2\theta\}\,d\theta
 =L\left(\frac12-\frac1\pi\right)<\frac L4.
\]
Here $\min\{\cos^2\theta,\sin^2\theta\}
=(1-|\cos2\theta|)/2$ and the mean of $|\cos2\theta|$ is $2/\pi$.
Consequently some direction and its associated deterministic mask satisfy
$c_{\mathbf I}^2\le L(1/2-1/\pi)<L/4$ by the Rayleigh quotient.
The random-mask lower bound is thus genuinely a high-probability
conditioning statement, not a disguised all-mask conditioning claim.
Each such mask still has positive probability $2^{-L}$ under the
specified observation law.

\paragraph{Noise normalization and limits of the inference.}
Let $y=P_{\mathbf I}A_Lx+e$ and let $\widehat x\in S^1$ be a global
least-squares minimizer using the exact map. On the certified event,
\[
 \|\widehat x-x\|\le
 \frac{2\|e\|}{\sqrt{L(1-\epsilon)/2}}.
\]
For a fixed total window noise budget $\|e\|\le E$, this decreases like
$L^{-1/2}$. But if each retained scalar error is bounded by $\sigma$, then
$\|e\|\le\sqrt L\sigma$, yielding only
$\|\widehat x-x\|\le2\sigma\sqrt{2/(1-\epsilon)}$, independent of $L$.
For the $\epsilon=1/2$ example this bound is $4\sigma$, not
$4\sigma/\sqrt L$. Equivalently, normalizing both map and data by
$\sqrt L$ gives full constant one and retained lower constant
$\sqrt{(1-\epsilon)/2}$, with no deterministic gain hidden in the stacking
convention. These bounds may be capped by the state-space diameter two.

A stochastic averaging gain requires extra assumptions. For illustration
only, if errors are independent of the mask, have conditional mean zero
and conditional covariance $\sigma^2I_L$, the exact-map, unconstrained
linear least-squares estimate $\widetilde x\in\mathbb R^2$ satisfies,
conditional on any mask in the certified event,
\[
 \mathbb E[\|\widetilde x-x\|^2\mid\mathbf I]
 =\sigma^2\operatorname{tr}[(A_L^\top P_{\mathbf I}^\top
 P_{\mathbf I}A_L)^{-1}]
 \le\frac{4\sigma^2}{L(1-\epsilon)}.
\]
This follows directly by computing the estimator covariance and bounding
both Gram eigenvalues below by $L(1-\epsilon)/2$. It concerns a known
linear observation model and an unconstrained estimator. Dynamical model error and
approximation or optimization residuals require additional terms.

\paragraph{Exact finite-window diagnostic.}
Set $\beta_t=(t-L+1)\omega$ and let $s_t=+1$ when the first coordinate is
retained and $s_t=-1$ when the second is retained. The retained Gram matrix
has eigenvalues
\[
 \lambda_\pm=\frac L2\pm\frac12
 \left|\sum_{t=0}^{L-1}s_t e^{2i\beta_t}\right|.
\]
Indeed, the first-coordinate row is $(\cos\beta_t,-\sin\beta_t)$,
whereas the second-coordinate row is $(\sin\beta_t,\cos\beta_t)$.
Their outer products equal $I_2/2$ plus or minus one half of the symmetric
matrix with entries $(\cos2\beta_t,-\sin2\beta_t;
-\sin2\beta_t,-\cos2\beta_t)$. Summing gives the displayed eigenvalues.
Since every direction occurs as a circle secant, $\sqrt{\lambda_-}$ is
the exact realized lower-secant constant. Numerical checks of this identity
are implementation diagnostics only; the probability certificate above
is analytic and does not rely on a Monte Carlo success rate.

\section{Invertible coordinates: proofs and limits of the construction}
\label{app:invertible_coordinates}

For $a\geq0$, define $h_a$ on each finite-dimensional Euclidean
coordinate space by
\begin{equation}
 [h_a(z)]_j=z_j+a\tanh(z_j).
 \label{eq:ha_definition}
\end{equation}
This coordinatewise example respects a known mask: evaluating a retained
coordinate does not require recovering any of the missing coordinates.
We will show that
\begin{equation}
 \|z-w\|_2\le\|h_a(z)-h_a(w)\|_2\le(1+a)\|z-w\|_2.
 \label{eq:ha_bilipschitz}
\end{equation}

\begin{proposition}[Global coordinatewise diffeomorphism]
For $a\geq0$, the map $h_a$ in Eq.~\eqref{eq:ha_definition} is a
$C^\infty$ diffeomorphism from $\mathbb R^D$ onto $\mathbb R^D$; the same conclusion holds for $h_a^{(q)}$ on $\mathbb R^q$.
Its Jacobian has singular values in $[1,1+a]$, its inverse is
$1$-Lipschitz, and Eq.~\eqref{eq:ha_bilipschitz} holds globally.
\end{proposition}
\begin{proof}
Write $s_a(u)=u+a\tanh u$. Differentiation gives
\[
 s_a'(u)=1+a\operatorname{sech}^2u,\qquad
 1\leq s_a'(u)\leq1+a.
\]
For $u<v$, the fundamental theorem of calculus gives
\[
 v-u\leq s_a(v)-s_a(u)
 =\int_u^v s_a'(r)\,dr\leq(1+a)(v-u).
\]
Thus $s_a$ is strictly increasing. Since $|a\tanh u|\leq a$,
$s_a(u)\to\pm\infty$ as $u\to\pm\infty$. The intermediate value theorem
therefore proves surjectivity, and strict monotonicity proves injectivity.
The inverse function theorem applies at every point because $s_a'$ never
vanishes. The resulting local smooth inverses agree, by uniqueness, to form
a global smooth inverse. Componentwise application proves the same for
$h_a$. Squaring the one-dimensional inequalities and summing over coordinates
yields Eq.~\eqref{eq:ha_bilipschitz}. Replacing $z,w$ by inverse images gives
\[
 \|h_a^{-1}(v)-h_a^{-1}(u)\|_2\leq\|v-u\|_2.
\]
Finally, $Dh_a(z)$ is diagonal with entries
$1+a\operatorname{sech}^2z_j$, which proves the singular-value claim.
For scalar inverse evaluation, $u=s_a(z)$ implies
$z\in[u-a,u+a]$; a monotone bracketed solve is therefore available.
If an approximate inverse $\widetilde z$ has scalar residual at most $\epsilon$,
then $|\widetilde z-s_a^{-1}(u)|\leq\epsilon$, by the lower Lipschitz bound.
\end{proof}

\begin{corollary}[Classical embedding-induced conjugacy]\label{cor:classical-conjugacy}
Let $M$ be compact and $F:M\to M$ continuous. If
$E:M\to\mathbb R^q$ is continuous and injective, then
$H=h_a^{(q)}\circ E$ is a homeomorphism from $M$ onto $S=H(M)$.
The map $G=H\circ F\circ H^{-1}:S\to S$ satisfies
\[
 G\circ H=H\circ F,\qquad G^n\circ H=H\circ F^n
 \quad(n\geq0).
\]
If $M$ is a smooth manifold and $E$ a $C^k$ embedding, then $H$ is a
$C^k$ diffeomorphism onto its image manifold. If, in addition, $F$ is a
$C^k$ diffeomorphism, so is $G$ on $S$. Here $q$ is the ambient
delay-coordinate dimension, which is distinct from $\dim M=d$.
\end{corollary}
\begin{proof}
The composition $H$ is continuous and injective. A continuous bijection
from a compact space to its image in a Hausdorff space has continuous inverse.
Indeed, a closed subset of $M$ is compact, its image is compact and hence
closed in $S$, so the inverse is continuous. For every $x\in M$,
\[
 G(H(x))=H(F(H^{-1}(H(x))))=H(F(x)).
\]
The iteration identity follows by induction: if it holds at $n$, then
$G^{n+1}H=G(HF^n)=(GH)F^n=HF^{n+1}$. In the smooth case,
composition with the ambient diffeomorphism $h_a^{(q)}$ preserves an embedding;
the claimed $C^k$ inverse is the inverse on $S$, not an inverse onto the
whole ambient space. The stated $C^k$ regularity and invertibility of $G$ then follow by
composition, with $G^{-1}=HF^{-1}H^{-1}$.
\end{proof}
This corollary applies only on the masks or success event for which the
embedding hypotheses have been established. It neither proves them anew nor
raises their success probability. If $F$ is not invertible, the conjugacy
identity still holds for continuous self-maps, but neither $F$ nor $G$ should
be described as invertible dynamics.

\paragraph{Continuous-time transformation.}
Suppose $z(t)$ solves $\dot z=f(z)$ and $y=h_a(z)$. The chain rule gives
\[
 \begin{aligned}
 \dot y&=Dh_a(z)\dot z=Dh_a(z)f(z)\\
 &=Dh_a(h_a^{-1}(y))f(h_a^{-1}(y))=:f_a(y).
\end{aligned}
\]
Consequently, wherever the flows exist,
$\Phi_a^t=h_a\Phi^t h_a^{-1}$. This exact formula requires the pushforward
field $f_a$; fitting an unrelated low-degree polynomial in $y$ does not
satisfy it identically. If $z$ denotes standardized observation coordinates,
this is a coordinate relation on that reconstructed domain and does not by
itself identify an unobserved physical state.

\paragraph{Why the discrete residual is not coordinate invariant.}
Consider the scalar field $f(z)=z$, step size one, and endpoints $z_0=1$,
$z_1=3$. The trapezoidal residual is
\[
 z_1-z_0-\tfrac12(f(z_0)+f(z_1))=2-2=0.
\]
Use the smooth diffeomorphism $h(z)=z+z^3$, whose derivative is
$1+3z^2>0$ and whose limits are $\pm\infty$. Its pushed-forward field is
$f_h(h(z))=(1+3z^2)z$. Hence $h(z_0)=2$, $h(z_1)=30$,
$f_h(h(z_0))=4$, and $f_h(h(z_1))=84$. The transformed residual is
\[
 30-2-\tfrac12(4+84)=28-44=-16\ne0.
\]
Thus exact equivalence of continuous flows does not imply exact equivalence
of their trapezoidal fitting objectives. The endpoints above describe a
trapezoidal step, not the exact time-one solution of $\dot z=z$.

\paragraph{Identity-preserving feature augmentation is a different construction.}
For a smooth $q:\mathbb R^d\to\mathbb R^s$, let
$\Psi(z)=(z,q(z))$. Projection onto the first $d$ components is an inverse
on $\Psi(\mathbb R^d)$, and
$D\Psi(z)v=(v,Dq(z)v)$ has trivial kernel. Therefore $\Psi$ is a smooth
embedding, a diffeomorphism onto its graph, not onto $\mathbb R^{d+s}$.
If $q$ is $L_q$-Lipschitz on a domain, then
\[
 \|\Psi(z)-\Psi(w)\|_2^2
 =\|z-w\|_2^2+\|q(z)-q(w)\|_2^2
\]
proves lower and upper Lipschitz constants $1$ and $\sqrt{1+L_q^2}$ there.
Retaining all linear coordinates while appending selected quadratic features
thus preserves information at the feature level. It does not prove a finite
polynomial closure of the induced dynamics, nor exact graph invariance of a
learned update. This elementary graph-embedding argument is not a novelty claim.

\paragraph{Time-dependent masks.}
If distinct successful masks yield coordinate maps $H_t$, the correct
coordinate relation is
\[
 G_t=H_{t+1}FH_t^{-1},\qquad G_tH_t=H_{t+1}F.
\]
A single autonomous conjugate $G$ is not obtained unless an additional
consistent identification has been established. For a mask sequence
$\omega$ with shift $\theta$, the analogous expression is
$G_\omega=H_{\theta\omega}FH_\omega^{-1}$. Fiberwise identities do not,
without joint continuity and a jointly continuous inverse, establish a
homeomorphism of an entire skew-product space. A fully erased finite window
with positive probability prevents an embedding guarantee for every mask
realization when the state space has more than one point.

\paragraph{Learned-model error through inverse coordinates.}
Suppose $\widehat G$ maps $S=H(M)$ into $S$, $H^{-1}$ is $L_{\rm inv}$-Lipschitz,
and
\[
 \sup_{x\in M}\|\widehat G(H(x))-H(F(x))\|\leq\epsilon.
\]
Define $\widetilde F=H^{-1}\widehat G H$. Then
\begin{align*}
 \|\widetilde F(x)-F(x)\|
 &=\|H^{-1}(\widehat G(H(x)))-H^{-1}(H(F(x)))\|\\
 &\leq L_{\rm inv}\|\widehat G(H(x))-H(F(x))\|\\
 &\leq L_{\rm inv}\epsilon.
\end{align*}
If $F$ is $L_F$-Lipschitz and the compared iterates stay in the declared
domain, their errors satisfy
\[
 e_{n+1}\leq L_Fe_n+L_{\rm inv}\epsilon,
 \qquad
 e_n\leq L_F^ne_0+L_{\rm inv}\epsilon\sum_{j=0}^{n-1}L_F^j.
\]
The second inequality follows by substituting the first repeatedly or by
induction. For the ambient map $h_a$ alone, $L_{\rm inv}=1$; for
$H=h_a\circ E$, the conditioning of $E^{-1}$ remains and cannot be dropped.
A homeomorphism alone need not have a Lipschitz inverse. Moreover, small
regression error is not exact topological conjugacy: the identity map on a
circle and an arbitrarily small nonzero rotation are uniformly close but
cannot be conjugate, because any conjugate of the identity is the identity.

\section{Whole-block erasures: a two-step torus control}
\label{app:dimension-torus}

Let $\mathbb S^1=\{u\in\mathbb C:|u|=1\}$, identifying each block with
two real coordinates. Compare the compact state spaces
\[
 M_2=\{(u_1,u_2,u_3)\in(\mathbb S^1)^3:u_3=u_1u_2\},
 \qquad M_3=(\mathbb S^1)^3
\]
in the common Euclidean chord metric inherited from $\mathbb C^3$.
Advance each block by $u_j\mapsto e^{i\omega_j}u_j$, with
$\omega_3=\omega_1+\omega_2$ modulo $2\pi$ for the graph model.
Each mask $B_t\subset\{1,2,3\}$ retains exactly two whole labeled blocks.
For a fixed one- or two-step mask sequence $B$, let $E_B$ be the
corresponding ordered observation map; define $c_B$ and $C_B$ as
the infimum and supremum of
$\|E_B(u)-E_B(v)\|/\|u-v\|$ over distinct states in the given space.
These are the exact lower and upper secant constants:
\[
\begin{array}{c|cc|cc}
 &\multicolumn{2}{c|}{M_2}&\multicolumn{2}{c}{M_3}\\
 \text{mask pattern}&c_B&C_B&c_B&C_B\\ \hline
 \text{one epoch} & 1/\sqrt{3} &1&0&1\\
 \text{two equal masks}&\sqrt{2/3}&\sqrt{2}&0&\sqrt{2}\\
 \text{two distinct masks}&1&\sqrt{5/3}&1&\sqrt{2}
\end{array}
\]
Consequently, every one- or two-step mask embeds $M_2$; a two-step
mask embeds $M_3$ exactly when the two masks differ. If the masks
are independently drawn uniformly from the three possible pairs,
six of the nine ordered pairs are distinct and the two-step success
probability on $M_3$ is $2/3$. Each successful fixed-mask map induces
the exact conjugate of the rotation on its image.

\paragraph{Proof of the secant claims.}
Write $d_j=|u_j-v_j|$, $S=\sum_{j=1}^3d_j^2>0$, and let $n_j$ be the
number of times block $j$ is retained. Every rotation is a blockwise isometry, so
\[
 \frac{\|E_B(u)-E_B(v)\|^2}{\|u-v\|^2}
 =\frac{\sum_{j=1}^3 n_jd_j^2}{S}.
\]
On $M_3$, any one block can vary with the others fixed, so the extrema
of the squared ratio are $\min_j n_j$ and $\max_j n_j$.
On $M_2$, if $\{m,a,b\}=\{1,2,3\}$, the product relation (or its
conjugated rearrangement) gives $d_m\le d_a+d_b$ by adding and
subtracting a mixed product and using unit modulus. Thus
$d_m^2\le 2(d_a^2+d_b^2)$, hence $0\le d_m^2/S\le2/3$.
The lower endpoint is attainable by keeping phase $m$ fixed and
varying the other two compatibly. For the upper endpoint, write
$(w_1,w_2,w_3)=(u_1,u_2,\overline{u_3})$, so $\prod_j w_j=1$.
Start at $(1,1,1)$ and perturb phases $a,b$ by $s$ and phase $m$ by
$-2s$. The ratio $d_m^2/S$ then tends to
$4/(4+1+1)=2/3$ as $s\to0$, since $|e^{ix}-1|/|x|\to1$.
The squared observation ratio is $1-d_m^2/S$ for a single mask
omitting $m$, $2(1-d_m^2/S)$ for a repeated mask omitting $m$, and
$1+d_m^2/S$ when the masks are distinct and share block $m$.
Taking the two endpoints proves the table. A positive lower
constant implies injectivity and a Lipschitz inverse on the image;
it also implies immersion by applying the bound to differentiable
curves through each point. The conjugate map is
$E_B\varphi E_B^{-1}$, where $\varphi$ is the rotation on the
respective state space. When the lower constant is zero on $M_3$
there is an unobserved free phase, which gives an actual collision.

\paragraph{A paired population diagnostic, not a dimension estimator.}
Draw independent starts from Haar measure on each torus, with masks
independent of the states. Every single-time retained pair has the
same product Haar law in both models: in $M_2$ this follows from the
unimodular changes of angle coordinates from $(\theta_1,\theta_2)$
to each retained pair. For runs with $B_0=\{1,2\}$ and $3\in B_1$,
an event of probability $2/9$ under the independent uniform mask law,
form the complex product
\[
 Z=u_1(0)u_2(0)\overline{u_3(1)}.
\]
On $M_2$ it is identically $e^{-i\omega_3}$; on $M_3$ it has mean zero
and $|Z|=1$. Given $n\ge1$ eligible independent runs, the graph mean
has modulus one, whereas on $M_3$ the empirical mean satisfies
$\mathbb E|\overline Z_n|^2=1/n$ by expanding the square and using
independence and centering. Thus
$\Pr(|\overline Z_n|\ge\delta)\le\min\{1,1/(n\delta^2)\}$ for
$\delta>0$. In particular one eligible sample does not distinguish
the models by modulus. Reusing windows from one trajectory does not
give independent starts, and the diagnostic can remain constant on
an $M_3$ trajectory when $\omega_1+\omega_2-\omega_3=0$ modulo
$2\pi$. The construction is a noiseless whole-block control with
labeled masks and a known model; it makes no claim for arbitrary
scalar-coordinate erasures, empirical dimension certification, or
a learned predictor.

\paragraph{Executable local-PCA control.}
The source archive includes a fixed-seed noiseless calculation with
\(8192\) independent phase starts per model and two exogenous uniform
masks per start. For dimension estimation on \emph{both} models we
keep only distinct-mask windows, so the algorithm can transport the
actually measured blocks back by the specified rotations and assemble
complete vectors without assuming either dimension or the graph
constraint. The equal-mask cases remain in the exact secant comparison.
There were \(5500\) complete windows per model. With radius \(0.65\)
and tail fraction \(0.04\), excluding the anchor itself and abstaining below 12 neighbors, Eq.~(1.1) of Ref.~\cite{dimLim2024}
returned \(2\) and \(3\), respectively, at all \(64\) tested anchors
per model. Minimum neighborhood counts were \(87\) and \(15\).
On \(1882\) eligible starts the separate phase-moment magnitudes were
\(1\) and \(0.01469\). These are results of one synthetic run;
the \(64\) local decisions reuse the same cloud and are not independent
trials. We do not assert that the numerical settings satisfy the
published estimator's finite-sample sufficient bounds. The calculation
tests neither inference of the rotation from incomplete observations
nor the dimension of an arbitrary partially observed attractor.

\end{document}